# On quantitative uniqueness for parabolic equations

Igor Kukavica and Quinn Le

ABSTRACT. We consider the quantitative uniqueness properties for a parabolic type equation $u_t - \Delta u = w(x,t)\nabla u + v(x,t)u$, when $v \in L_t^{p_2} L_x^{p_1}$ and $w \in L_t^{q_2} L_x^{q_1}$, with a suitable range for exponents $p_1$, $p_2$, $q_1$, and $q_2$. We prove a strong unique continuation property and provide a pointwise in time observability estimate.          July 27, 2021

## CONTENTS

## 1. Introduction

We address quantitative uniqueness for solutions of the parabolic equation

$$u_t - \Delta u = w(x,t)\nabla u + v(x,t)u, \qquad (x,t) \in \Omega \times [T_0, T_0 + T],$$

where $v$ and $w$ belong to the mixed Lebesgue spaces $v \in L_t^{p_2} L_x^{p_1}$ and $w \in L_t^{q_2} L_x^{q_1}$, respectively. In particular, we obtain bounds on the order of vanishing which are algebraic in the corresponding norms of $v$ and $w$. Additionally, we obtain pointwise in time doubling estimates, i.e., inequalities of the form

$$\|u(\cdot,t)\|_{L^2(\Omega)} \leq M_\delta \|u(\cdot,t)\|_{L^2(B_\delta(x_0))}$$

for an arbitrary $\delta > 0$ under the Lebesgue conditions on $v$ and $w$.

The unique continuation for PDEs has a rich history (cf. the review papers by Kenig [K1, K2] and Vessella [V]), so we only mention several results pertaining to this paper. In [JK], Jerison and Kenig proved that the second order elliptic equation has the strong unique continuation property (i.e., is identically zero if it vanishes to an infinite order at a point) if $w = 0$ and $v \in L^{n/2}$, with a sufficiently small $L^{n/2}$ norm, which is a sharp result. The parabolic counterpart was obtained by Escauriaza and Vega [EV] (see also [E] for a previous unique continuation result when $v \in L^\infty L^p$ with $p > n/2$). The difficult case when $w$ is nonzero was addressed by Koch and Tataru for the elliptic case in [KT1] and the parabolic case in [KT2]. In particular, they obtained the strong unique continuation for $v \in L_t^1 L_x^\infty + L_t^\infty L_x^{n/2}$ with the norm sufficiently small and $w \in L_{x,t}^{n+2}$. All the mentioned works rely on suitable Carleman type estimates and lead to observability type estimates on the space-time rectangles. For other works on the frequency approach to the unique continuation, see [A, AN, Ch, Ku1, Ku2], for the works related to Dirichlet quotients, see [A, CFNT, FS] and for other related works, cf. [AE, AMRV, An, AV, BC, B, BK, CRV, DF1, DF2, DZ1, DZ2, EF, EFV, EVe, H1, H2, K3, KSW, L, M, Z].

In this paper (cf. Theorems (2.1)–(2.3) below), we obtain explicit algebraic observability estimates for fixed time (i.e., not only on space-time rectangles) under the assumptions on the coefficients $v \in L^\infty L^p$ and $w \in L^\infty L^q$, where $p > 2n/3$ and $q > 2n$. More general conditions $v \in L_t^{p_2} L_x^{p_1}$ and $w \in L_t^{q_2} L_x^{q_1}$ under certain range are addressed in Theorems 2.1–2.3. While the results cited in the second paragraph use the Carleman estimates, we rely on the frequency function approach, developed in [**Al**, **GL**] for the elliptic and [**Kur**, **P**] for the parabolic equations. The main idea of this approach for parabolic equations is the logarithmic convexity of the frequency function

$$Q(t) = \frac{|t| \int_{\mathbb{R}^n} |\nabla u(x,t)|^2 G(x,t)\, dx}{\int_{\mathbb{R}^n} u(x,t)^2 G(x,t)\, dx} \tag{1.1}$$

for the heat equation. In (1.1), $G$ is the $(4\pi)^{n/2}$-multiple of the backward Gaussian kernel, i.e.,

$$G(x,t) = \frac{1}{|t|^{n/2}} e^{|x|^2/4t}.$$

Another reformulation of the idea is to use the similarity variables (cf. (3.8) and (3.9) below) and obtain a logarithmic convexity of the unweighted norm [**C**, **Ku3**]. In [**CK**], the approach was used to obtain an estimate for an order of vanishing $C(\|w\|_{L^\infty}^2 + \|v\|_{L^\infty}^{2/3})$, which is, at least for complex valued coefficients, sharp [**CKW1**, **CKW2**]. Poon and Kurata showed that the frequency approach leads to unique continuation for $p > n$ and $q = \infty$.

In this paper, we show that this may be improved to $p > 2n/3$ and $q > 2n$. The improved range is obtained by three main devices. The first is to find the point in space where the frequency function is the smallest and translate the equation so it starts at that point (cf. Lemma 3.1 below); this idea has been introduced in [**Ku4**, **CK**]. The second device is to use the embedding theorems with Laplacian and use (3.24)–(3.25) below to bound the parts containing $v$ and $w$. The third device is to use the finiteness of the integral in (3.47) below, which then allows us to show the convergence of the quantity under the integral. Note that we obtain an explicit algebraic bound on the order of vanishing, which is a constant multiple of

$$\|v\|_{L_t^\infty L_x^p}^{2/(3-2n/p)} + \|w\|_{L_t^\infty L_x^q}^{2/(1-2n/q)}.$$

When setting $p = q = \infty$, the estimate reduces to the sharp bound from [**CK**]. We also provide a pointwise estimate in time for a better understanding of the behavior of solution $u$. In particular, for all $\delta_0 \in (0,1]$, we have

$$\|u(\cdot,t)\|_{L^2(\mathbb{R}^n)} \lesssim e^P \|u(\cdot,t)\|_{L^2(B_{\delta_0}(0))} \tag{1.2}$$

for all $t \in (T_0, T_0 + T)$, where $P$ is a polynomial depending only on $n$, $q$, $p$, $\delta_0$, $\|v\|_{L^p}$ and $\|w\|_{L^q}$. The explicit formula for $P$ can be found in Lemma 4.2 below. Note that the estimates of the type (1.2) are an essential ingredient when considering qualitative properties of solutions of evolutionary parabolic PDE. For instance, they are needed when considering the size of the zero set of solution at time $t$ or, more generally, complexity of a graph of a function at a time $t$ ([**Ku3**, **Ku4**]). Note finally that compared to [**CK**], we reduce the necessary regularity for the solution $u$ to a simple boundedness, rather than differentiability.

The paper is structured as follows. In Section 2, we provide the setup of the problem and state the two main theorems on the order of vanishing with periodic boundary conditions and the analog for $\mathbb{R}^n$. For simplicity, the theorems address periodic boundary conditions. The case of $\mathbb{R}^n$ under a natural non-growth condition on the solution is considered in Theorem 2.3 and proven in Section 4. Theorems 2.1–2.3 are concerned with the case when the Lebesgue exponent of $t$ is infinity, while the general range is addressed in Theorem 6.1 in Section 6.

## 2. The main result

We address the quantitative uniqueness of a nontrivial solution $u \in L^\infty(I, L^2(\mathbb{T}^n)) \cap L^2(I, H^1(\mathbb{T}^n))$ of the equation

$$\begin{aligned} u_t - \Delta u &= w(x,t) \nabla u + v(x,t) u \\ u(x,T_0) &= u_0(t) \end{aligned} \tag{2.1}$$

defined for $(x,t) \in \mathbb{R}^n \times I$ where $I = [T_0, T_0 + T]$ is a given time interval and $n \geq 3$. The theorems are also valid for $n = 1, 2$ with minor changes; cf. Remark 2.4 below. We assume that $u$, $v$, and $w$ are $\Omega$-periodic, where $\Omega = \mathbb{T}^n = [-1/2, 1/2]^n$ and that they satisfy

$$\|v(\cdot, t)\|_{L^p(\Omega)} \leq M_0 \tag{2.2}$$

and

$$\|w(\cdot, t)\|_{L^q(\Omega)} \leq M_1, \tag{2.3}$$

for all $t \in I$. For simplicity, and since we are interested in the regime of $M_1$ and $M_2$ large, we assume $M_1, M_2 \geq 1$. Let $O_{(x_0,t_0)}(u)$ be the vanishing order of $u$ at $(x_0, t_0)$, which is defined as the largest integer $d$ such that

$$\|u\|_{L^2(Q_r(x_0,t_0))} = O(r^{d+(n+2)/2}) \quad \text{as } r \to 0^+$$

where

$$Q_r(x_0, t_0) = \{(x,t) \in \mathbb{R}^n \times \mathbb{R} : |x - x_0| < r, -r^2 < t - t_0 < 0\} \tag{2.4}$$

stands for the parabolic cylinder centered at $(x_0, t_0)$ with the radius $r$. For $t \in I$, denote by

$$q_D(t) = \frac{\|\nabla u(\cdot, t)\|_{L^2(\Omega)}^2}{\|u(\cdot, t)\|_{L^2(\Omega)}^2} \tag{2.5}$$

the Dirichlet quotient of $u$ at time $t$. We assume that $\|u(\cdot, t)\|_{L^2(\Omega)}$ is nonzero for all $t \in I$. We also assume

$$q_0 = \sup_{t \in I} q_D(t) < \infty. \tag{2.6}$$

The following is the main result of this paper; cf. also Theorem 6.1 for $v \in L_t^{p_2} L_x^{p_1}(\Omega \times I)$ and $w \in L_t^{q_2} L_x^{q_1}(\Omega \times I)$.

THEOREM 2.1. *Let $u \in L_{x,t}^\infty(\mathbb{T}^n \times I)$ be a solution of* (2.1) *with $v$ and $w$ satisfying* (2.2) *and* (2.3) *such that*

$$p > \frac{2n}{3}$$

*and*

$$q > 2n.$$

*Then for all $(x_0, t_0) \in \mathbb{T}^n \times I$, the vanishing order of $u$ at $(x_0, t_0)$ satisfies*

$$O_{(x_0,t_0)}(u) \lesssim M_0^a + M_1^b, \tag{2.7}$$

*where $a = 2/(3 - 2n/p)$ and $b = 2/(1 - 2n/q)$, with the implicit constant in* (2.7) *depending on $q_0$ and $I$.*

We emphasize that the dimension is considered fixed, so we allow all constants to depend on $n$ without mention. In the next statement, we provide a pointwise in time doubling property of $u$ under the conditions of Theorem 2.1.

THEOREM 2.2. *Under the conditions of Theorem 2.1, we have*

$$\|u(\cdot, t)\|_{L^2(\Omega)} \lesssim e^{P(n, \delta_0, M_0, M_1, p, q)} \|u(\cdot, t)\|_{L^2(B(0, \delta_0))},$$

*for all $t \in (T_0, T_0 + T)$, where $P$ is a polynomial.*

Now, consider $u$ to be given on $\mathbb{R}^n$, instead of $\mathbb{T}^n$. Suppose that it satisfies the condition

$$\int_{\mathbb{R}^n} u(x, t)^2 \, dx \leq K \int_{B(0,1)} u(x, t)^2 \, dx, \qquad t \in [T_0, T_0 + T], \tag{2.8}$$

for some constant $K$. In this case, we obtain the following analogue of Theorems 2.1 and 2.2.

THEOREM 2.3. *Let $u \in L_t^\infty L_x^\infty$ be a solution of* (2.1) *such that the condition* (2.8) *holds with $v$ and $w$ satisfying* (2.2) *and* (2.3). *Assume additionally that*

$$q > 2n$$

*and*

$$p > \frac{2n}{3}.$$

*Then for all $(x_0, t_0) \in \mathbb{R}^n \times I$, the vanishing order of $u$ at $(x_0, t_0)$ satisfies*

$$O_{(x_0,t_0)}(u) \lesssim M_0^a + M_1^b, \tag{2.9}$$

*where $a = 2/(3 - 2n/p)$ and $b = 2/(1 - 2n/q)$, with the implicit constant in* (2.9) *depending on $q_0$, $K$, and $I$. Moreover, for $\delta_0 \in (0, 1]$, we have*

$$\|u(\cdot, t)\|_{L^2(\Omega)} \lesssim e^{P(n, \delta_0, M_0, M_1, q, p)} \|u(\cdot, t)\|_{L^2(B(0, \delta_0))}, \tag{2.10}$$

*for all $t \in (T_0, T_0 + T)$, where $P$ is a polynomial.*

The theorem is proven in Section 5 below.

REMARK 2.4. In the theorems above, we assumed $n \geq 3$. In this remark, we address the changes when $n = 1$ or $n = 2$. If $n = 2$, the results are exactly the same with small modifications in the proof; the modifications are pointed out in the proof of Theorem 2.1. For the case $n = 1$, we need to assume additionally that $p \geq 2$ and $q \geq 4$.

## 3. Proof on the statement for quantitative uniqueness

This section is devoted to the proof of Theorem 2.1. We first start with the case where $v$, $w$, and $u_0$ are smooth and then use approximation argument to prove the theorem under the boundedness assumption. Thus assume for now that $v$, $w$, and $u_0$ are smooth. By translation and rescaling, we may assume that $I = [-1, 0]$ and $(x_0, t_0) = (0, 0)$. Suppose additionally that $u$ is periodic in $x$ with respect to $\Omega = [-1/2, 1/2]^n$. For any starting time $-\epsilon$, we claim that we can choose a starting point $x_\epsilon$ that is controlled by $q_0$.

LEMMA 3.1. *Let $u$ be a solution of* (2.1) *that is periodic in $x$ with respect to $\Omega = [-1/2, 1/2]^n$ where $v$ and $w$ satisfying* (2.2) *and* (2.3). *For any $\epsilon \in (0, T)$, there exist $x_\epsilon \in \Omega$*

$$\frac{\epsilon \int_{\mathbb{R}^n} |\Delta u(x_\epsilon + y, -\epsilon)|^2 G(y, -\epsilon)\, dy}{\int_{\mathbb{R}^n} u(x_\epsilon + y, -\epsilon)^2 G(y, -\epsilon)\, dy} \leq 2\epsilon q_D(-\epsilon),$$

*where*

$$G(x, t) = \frac{1}{|t|^{n/2}} e^{|x|^2/4t}, \qquad x \in \mathbb{R}^n, \quad t < 0. \tag{3.1}$$

The lemma is proven in [**CK, Ku4**]; we provide a short argument for the sake of completeness. In the sequel, we reserve $\epsilon$ for the time provided in Lemma 3.1. We always tacitly assume that $\epsilon \in (0, T)$.

PROOF OF LEMMA 3.1. Assume that we have

$$2q_D(-\epsilon) \int_{\mathbb{R}^n} u(x+y, -\epsilon)^2 G(y, -\epsilon)\, dy \leq \int_{\mathbb{R}^n} |\nabla u(x+y, -\epsilon)|^2 G(y, -\epsilon)\, dy, \tag{3.2}$$

for all $x \in \Omega$. With a simple change of variables, (3.2) becomes

$$\int_{\mathbb{R}^n} u(y, -\epsilon) G(y - x, -\epsilon)\, dy \leq \frac{1}{2q_D(-\epsilon)} \int_{\mathbb{R}^n} |\nabla u(y, -\epsilon)|^2 G(y - x, -\epsilon)\, dy, \qquad x \in \Omega. \tag{3.3}$$

Integrating (3.3) over $\Omega$, we obtain

$$\begin{aligned}
\int_\Omega \int_{\mathbb{R}^n} u(y, -\epsilon)^2 G(y, -\epsilon)\, dy\, dx &= \int_\Omega \sum_{j \in \mathbb{Z}^n} \int_{j+\Omega} u(y, -\epsilon)^2 G(y - x, -\epsilon)\, dy\, dx \\
&= \int_\Omega \sum_{j \in \mathbb{Z}^n} \int_\Omega u(y', -\epsilon)^2 G(y' + j - x, -\epsilon)\, dy'\, dx \\
&= \int_\Omega u(y', -\epsilon)^2 \left( \sum_{j \in \mathbb{Z}^n} \int_\Omega G(y' + j - x, -\epsilon)\, dx \right) dy' \\
&= \int_\Omega u(y', -\epsilon)^2 \int_{\mathbb{R}^n} G(y' - x, -\epsilon)\, dx\, dy' = (4\pi)^{n/2} \int_\Omega u(y', -\epsilon)^2\, dy',
\end{aligned} \tag{3.4}$$

where we used $\int_{\mathbb{R}^n} G(y - x, -\epsilon)\, dx = (4\pi)^{n/2}$ in the last equality. Similarly, we have

$$\int_\Omega \int_{\mathbb{R}^n} |\nabla u(y, -\epsilon)|^2 G(y - x, -\epsilon)\, dy\, dx = (4\pi)^{n/2} \int_\Omega |\nabla u(y', -\epsilon)|^2\, dy'. \tag{3.5}$$

By (3.3), (3.4), and (3.5), we obtain

$$2q_D(-\epsilon) \int_\Omega u(y, -\epsilon)^2\, dy \leq \int_\Omega |\nabla u(y, -\epsilon)|^2\, dy = q_D(-\epsilon) \int_\Omega u(y, -\epsilon)^2\, dy.$$

We have a contradiction since $\int_\Omega u(y, -\epsilon)^2\, dy > 0$. Therefore, (3.2) cannot hold for all $x \in \Omega$. Thus Lemma 3.1 follows. $\square$

We proceed with a change of variables so that at the time $-\epsilon$, given by Lemma 3.1, the solution $u$ starts at $x_\epsilon$ with $x_\epsilon$ also as in Lemma 3.1. In particular, we let

$$u(x, t) = \bar{u}\left(x + \frac{x_\epsilon}{\epsilon} t, t\right).$$

It is not difficult to check that $\bar{u}$ solves the equation

$$\partial_t \bar{u} - \Delta \bar{u} = -\frac{1}{\epsilon} x_\epsilon \nabla \bar{u} + w \nabla \bar{u} + v \bar{u}. \tag{3.6}$$

Denoting $a = -x_\epsilon/\epsilon$ and writing $u$ instead of $\bar{u}$, the equation (3.6) becomes

$$\partial_t u - \Delta u = a_j \partial_j u + w_j \partial_j u + vu. \tag{3.7}$$

We now proceed with a change of variable

$$u(x,t) = e^{|x|^2/8(-t)} U\left(\frac{x}{\sqrt{-t}}, -\log(-t)\right), \tag{3.8}$$

that is,

$$U(y,\tau) = e^{-|y|^2/8} u(ye^{-\tau/2}, -e^{-\tau}) \tag{3.9}$$

with $\tau = -\log(-t)$. Also, let

$$V(y,\tau) = v(ye^{-\tau/2}, -e^{-\tau}), \qquad (y,\tau) \in \mathbb{R}^n \times [\tau_0, \infty)$$

and

$$W_j(y,\tau) = w_j(ye^{-\tau/2}, -e^{-\tau}), \qquad (y,\tau) \in \mathbb{R}^n \times [\tau_0, \infty),$$

for $j = 1, ..., n$, where

$$\tau_0 = \log(1/\epsilon). \tag{3.10}$$

Then

$$\partial_\tau U + HU = e^{-\tau/2}(a_j y_j U + a_j \partial_j U) + e^{-\tau/2}(y_j W_j(y,\tau) U + W_j(y,\tau) \partial_j U) + e^{-\tau} V(y,\tau) U, \tag{3.11}$$

where

$$HU = -\Delta U + \left(\frac{|y|^2}{16} - \frac{n}{4}\right) U$$

with the initial data

$$U(y,\tau_0) = U\left(y, \log\frac{1}{\epsilon}\right) = e^{-|y|^2/8} u\left(\frac{y}{\sqrt{\epsilon}}, -\epsilon\right).$$

Now, define the frequency function

$$Q(\tau) = \frac{(HU, U)_{L^2(\mathbb{R}^n)}}{\|U\|^2},$$

where we denote

$$\|\cdot\| = \|\cdot\|_{L^2(\mathbb{R}^n)} \tag{3.12}$$

throughout. Letting

$$A(\tau)U = HU - e^{-\tau/2} a_j y_j U \tag{3.13}$$

and

$$\bar{Q}(\tau) = \frac{(A(\tau)U, U)_{L^2(\mathbb{R}^n)}}{\|U\|^2} = Q(\tau) - \frac{e^{-\tau/2} a_j}{\|U\|^2} \int_{\mathbb{R}^n} y_j U^2 \, dy, \tag{3.14}$$

we may rewrite (3.11) as

$$\partial_\tau U + (A(\tau) - \bar{Q}(\tau)I)U + \bar{Q}(\tau)U = e^{-\tau/2} a_j \partial_j U + e^{-\tau/2}(y_j W_j U + W_j \partial_j U) + e^{-\tau} VU. \tag{3.15}$$

For simplicity of notation, denote

$$\tilde{U} = \frac{U}{\|U\|},$$

so that $\|\tilde{U}\| = 1$. We now show that the modified frequency function $\bar{Q}$ is bounded.

LEMMA 3.2. *Under the assumptions of Theorem 2.1, for all $\tau \geq \tau_0$, the modified frequency function satisfies*

$$\bar{Q}(\tau) \lesssim M_0^a + M_1^b, \tag{3.16}$$

*where $a = 2/(3 - 2n/p)$ and $b = 2/(1 - 2n/q)$.*

PROOF OF LEMMA 3.2. Without loss of generality, assume that $M_0, M_1 \geq 1$. Dividing (3.15) by $\|U\|$ and then taking the inner product with $(A(\tau) - \bar{Q}(\tau)I)\tilde{U}$, we obtain, using $((A(\tau) - \bar{Q}(\tau)I)\tilde{U}, \tilde{U}) = 0$,

$$\frac{1}{\|U\|}\bigl(\partial_\tau \tilde{U}, (A(\tau) - \bar{Q}(\tau)I)\tilde{U}\bigr) + \bigl\|(A(\tau) - \bar{Q}(\tau)I)\tilde{U}\bigr\|^2 = \bigl(F(\tilde{U}), (A(\tau) - \bar{Q}(\tau)I)\tilde{U}\bigr), \qquad (3.17)$$

where $F(\tilde{U})$ denotes the right-hand side of (3.15) divided by $\|U\|$, i.e.,

$$F(\tilde{U}) = e^{-\tau/2} a_j \partial_j \tilde{U} + e^{-\tau/2}(y_j W_j \tilde{U} + W_j \partial_j \tilde{U}) + e^{-\tau} V \tilde{U}.$$

From (3.17), we get

$$\frac{1}{\|U\|}\bigl(\partial_\tau U, (A(\tau) - \bar{Q}(\tau)I)\tilde{U}\bigr) + \bigl\|(A(\tau) - \bar{Q}(\tau)I)\tilde{U}\bigr\|^2 = \bigl(F(\tilde{U}), (A(\tau) - \bar{Q}(\tau)I)\tilde{U}\bigr). \qquad (3.18)$$

By differentiating (3.14), we have

$$\bar{Q}'(\tau) = 2(\partial_\tau \tilde{U}, A(\tau)\tilde{U}) + (\tilde{U}, A'(\tau)\tilde{U}) - 2(\partial_\tau \tilde{U}, \tilde{U})(A(\tau)\tilde{U}, \tilde{U})$$

from where we use that $A(\tau)$ is a symmetric operator with the derivative $A'(\tau)U = \frac{1}{2}e^{-\tau/2} a_j y_j U$ to obtain

$$\begin{aligned}
\frac{1}{2}\bar{Q}'(\tau) &= (\partial_\tau \tilde{U}, A(\tau)\tilde{U}) + \frac{1}{4}e^{-\tau/2}(\tilde{U}, a_j y_j \tilde{U}) - (\partial_\tau \tilde{U}, \tilde{U})(A(\tau)\tilde{U}, \tilde{U}) \\
&= (\partial_\tau \tilde{U}, (A(\tau) - \bar{Q}(\tau)I)\tilde{U}) + \frac{1}{4}e^{-\tau/2}(\tilde{U}, a_j y_j \tilde{U}).
\end{aligned} \qquad (3.19)$$

Substituting (3.18) into (3.19), we get

$$\begin{aligned}
\frac{1}{2}\bar{Q}'(\tau) + \|(A(\tau) - \bar{Q}(\tau)I)\tilde{U}\|^2 &= \frac{1}{4}e^{-\tau/2}(a_j y_j \tilde{U}, \tilde{U}) + e^{-\tau/2} a_j \bigl(\partial_j \tilde{U}, (A(\tau) - \bar{Q}(\tau)I)\tilde{U}\bigr) \\
&\quad + \bigl(e^{-\tau/2}(y_j W_j \tilde{U} + W_j \partial_j \tilde{U}) + e^{-\tau} V \tilde{U}, (A(\tau) - \bar{Q}(\tau)I)\tilde{U}\bigr).
\end{aligned} \qquad (3.20)$$

Note that $\int \tilde{U}\partial_j U = 0$ and $\int \Delta U \partial_j U = 0$ since $v$, $w$, and $u_0$ are assumed smooth. Integrating by parts the second term on the right-hand side of (3.20), we obtain

$$\begin{aligned}
\frac{1}{2}\bar{Q}'(\tau) &= \frac{1}{4}e^{-\tau/2} a_j \int y_j \tilde{U}^2 \, dy - \frac{1}{16}e^{-\tau/2} a_j \int y_j \tilde{U}^2 \, dy + \frac{1}{2}e^{-\tau}|a|^2 \\
&\quad + e^{-\tau/2} \int (y_j W_j \tilde{U} + W_j \partial_j \tilde{U} + e^{-\tau/2} V \tilde{U})(A(\tau) - \bar{Q}(\tau)I)\tilde{U} \, dy \\
&= I_1 + I_2 + I_3 + I_4.
\end{aligned}$$

Using the Cauchy-Schwarz inequality and $\|\tilde{U}\| = 1$, we have

$$I_1 \leq e^{-\tau/2}|a| \int |y|\tilde{U}^2 \, dy \leq e^{-\tau/2}|a| \left(\int |y|^2 \tilde{U}^2 \, dy\right)^{1/2},$$

and similarly

$$I_2 + I_3 \leq e^{-\tau/2}|a| \left(\int |y|^2 \tilde{U}^2 \, dy\right)^{1/2} + e^{-\tau}|a|^2.$$

In order to estimate $I_4$, we first claim

$$\|D^2 U\|^2 \lesssim \|HU\|^2 + \|U\|^2. \qquad (3.21)$$

To prove (3.21), we expand $\|HU\|^2$ as

$$\begin{aligned}
\|HU\|^2 &= \int_{\mathbb{R}^n} \left(-\Delta U + \left(\frac{|y|^2}{16} - \frac{n}{4}\right) U\right)^2 dy \\
&= \int_{\mathbb{R}^n} (\Delta U)^2 + \left(\frac{|y|^2}{16} - \frac{n}{4}\right)^2 U^2 \, dy - 2\int_{\mathbb{R}^n} \left(\frac{|y|^2}{16} - \frac{n}{4}\right) U \Delta U \, dy \\
&= \int_{\mathbb{R}^n} \left((\Delta U)^2 + \left(\frac{|y|^2}{16} - \frac{n}{4}\right)^2 U^2\right) dy + 2\int_{\mathbb{R}^n} \partial_j U \partial_j U \left(\frac{|y|^2}{16} - \frac{n}{4}\right) dy + \frac{1}{4} \int_{\mathbb{R}^n} y_j U \partial_j U \, dy,
\end{aligned}$$

from where

$$\|HU\|^2 = \int_{\mathbb{R}^n} \left( (\Delta U)^2 + \left(\frac{|y|^2}{16} - \frac{n}{4}\right)^2 U^2 \right) dy + 2\int_{\mathbb{R}^n} |\nabla U|^2 \left(\frac{|y|^2}{16} - \frac{n}{4}\right) dy - \frac{n}{8} \int_{\mathbb{R}^n} U^2 \, dy$$

$$\geq \int_{\mathbb{R}^n} \left( (\Delta U)^2 + \left(\frac{|y|^2}{16} - \frac{n}{4}\right)^2 U^2 \right) dy - \frac{n}{2} \int |\nabla U|^2 \, dy - \frac{n}{8} \int_{\mathbb{R}^n} U^2 \, dy.$$

Using $n/4 - |y|^2/16 \leq n/4$ in the third integral, we get

$$\|\Delta U\|^2 \lesssim \|HU\|^2 + \|U\|^2 + \|\nabla U\|^2.$$

By Sobolev's and the Cauchy-Schwarz inequalities, we get

$$\|\nabla U\|^2 \leq \|U\|\|\Delta U\| \leq \frac{\|U\|^2 + \|\Delta U\|^2}{2},$$

and then using $\|D^2 U\| \lesssim \|\Delta U\|$, we obtain (3.21).

We now estimate $\|V\tilde{U}\|$. By Hölder's and Sobolev's inequalities, we have

$$\|V\tilde{U}\| \lesssim \|V\|_{L^p(\mathbb{R}^n)} \|\tilde{U}\|_{L^{2p/(p-2)}(\mathbb{R}^n)} \lesssim \|V\|_{L^p(\mathbb{R}^n)} \|\tilde{U}\|^{1-\alpha} \|D^2 \tilde{U}\|^{\alpha} \tag{3.22}$$

with

$$\alpha = \frac{n}{2p}. \tag{3.23}$$

Note that we need $p > 2$ in (3.22). This automatically holds if $n \geq 3$ and $p > 2n/3$. If $n = 2$, $p$ can be 2 since we only require $p > 2n/3 = 4/3$. In this case, we can simply replace $2p/(p-2)$ by $\infty$ and use Sobolev's inequalities for $L^{\infty}$. This estimate also works for the case $n = 1$ since we have assumed additionally that $p \geq 2$ if $n = 1$. By (3.21), we have

$$\|V\tilde{U}\| \lesssim \|V\|_{L^p(\mathbb{R}^n)} \|\tilde{U}\|^{1-\alpha} (\|H\tilde{U}\|^2 + \|\tilde{U}\|^2)^{\alpha/2} \lesssim \|V\|_{L^p(\mathbb{R}^n)} (\|H\tilde{U}\|^{\alpha} + 1). \tag{3.24}$$

In order to estimate the last expression, we write

$$\|H\tilde{U}\| \leq \|H\tilde{U} - A(\tau)\tilde{U}\| + \|(A(\tau) - \bar{Q}(\tau)I)\tilde{U}\| + \|\bar{Q}(\tau)\tilde{U}\|$$
$$\lesssim e^{-\tau/2}|a|\|y\tilde{U}\| + \|(A(\tau) - \bar{Q}(\tau)I)\tilde{U}\| + \bar{Q}(\tau), \tag{3.25}$$

where we used (3.13) in the last step. Since

$$\|\nabla \tilde{U}\| + \|y\tilde{U}\| \lesssim \bar{Q}(\tau)_+^{1/2} + e^{-\tau/2}|a| + 1, \tag{3.26}$$

([Ku4, p.780]), as one may readily check, we get

$$\|H\tilde{U}\| \lesssim e^{-\tau/2}|a|(\bar{Q}(\tau)_+^{1/2} + e^{-\tau/2}|a| + 1) + \|(A(\tau) - \bar{Q}(\tau)I)\tilde{U}\| + |\bar{Q}(\tau)|. \tag{3.27}$$

Since $\|V\|_{L^p(\mathbb{R}^n)} \leq M_0 e^{\alpha\tau}$ by a simple change of variables, we obtain from (3.24) and (3.27) the inequality

$$\|V\tilde{U}\| \lesssim M_0 e^{\alpha\tau} \left( e^{-\alpha\tau/2}|a|^{\alpha} \left( \bar{Q}(\tau)_+^{\alpha/2} + e^{-\alpha\tau/2}|a|^{\alpha} + 1 \right) + \|(A(\tau) - \bar{Q}(\tau)I)\tilde{U}\|^{\alpha} + |\bar{Q}(\tau)|^{\alpha} + 1 \right). \tag{3.28}$$

Next, to estimate $\|y_j W_j \tilde{U} + W_j \partial_j \tilde{U}\|$, we first bound $\|W_j \partial_j \tilde{U}\|$ and then $\|y_j W_j \tilde{U}\|$. By Hölder's and Sobolev's inequalities again, we have, also using $\|W\|_{L^p(\mathbb{R}^n)} \lesssim M_1 e^{n\tau/2q} = M_1 e^{(\beta - 1/2)\tau}$,

$$\|W_j \partial_j \tilde{U}\| \lesssim \|W\|_{L^q(\mathbb{R}^n)} \|\nabla \tilde{U}\|_{L^{2q/(q-2)}(\mathbb{R}^n)} \lesssim M_1 e^{(\beta - 1/2)\tau} \|\tilde{U}\|^{1-\beta} \|D^2 \tilde{U}\|^{\beta}$$

with

$$\beta = \frac{n}{2q} + \frac{1}{2}. \tag{3.29}$$

Similarly as in (3.25), we have

$$\|W_j \partial_j \tilde{U}\| \lesssim M_1 e^{(\beta - 1/2)\tau} \left( \|H\tilde{U}\|^{\beta} + 1 \right)$$
$$\lesssim M_1 e^{(\beta - 1/2)\tau} \left( e^{-\beta\tau/2}|a|^{\beta}(\bar{Q}(\tau)_+^{\beta/2} + e^{-\beta\tau/2}|a|^{\beta} + 1) + \|(A(\tau) - \bar{Q}(\tau)I)\tilde{U}\|^{\beta} + \bar{Q}(\tau)^{\beta} + 1 \right). \tag{3.30}$$

Next,

$$\begin{aligned}\|y_j W_j \tilde{U}\| &\lesssim \|W\|_{L^q(\mathbb{R}^n)} \|y\tilde{U}\|_{L^{2q/(q-2)}(\mathbb{R}^n)} \lesssim M_1 e^{(\beta-1/2)\tau} \|y\tilde{U}^{1/2}\|_{L^4(\mathbb{R}^n)} \|\tilde{U}^{1/2}\|_{L^{4q/(q-4)}(\mathbb{R}^n)}\\ &= M_1 e^{(\beta-1/2)\tau} \|y^2 \tilde{U}\|^{1/2} \|\tilde{U}\|_{L^{2q/(q-4)}(\mathbb{R}^n)}^{1/2} \lesssim M_1 e^{(\beta-1/2)\tau} \|y^2 \tilde{U}\|^{1/2} \|\Delta \tilde{U}\|^{n/2q} \|\tilde{U}\|^{1/2-n/2q} \\ &\lesssim M_1 e^{(\beta-1/2)\tau} \|y^2 \tilde{U}\|^{1/2} \|\Delta \tilde{U}\|^{n/2q},\end{aligned} \quad (3.31)$$

where we used (3.29) in the last inequality. Similarly to estimate $\|V\tilde{U}\|$, we need $q \geq 4$ in (3.31). This is obviously true if $n \geq 2$ and $q > 2n$. If $n = 1$, (3.31) still holds since we assumed $q \geq 4$. Recall that $HU = -\Delta U + (|y|^2/16 - n/4)U$, from where $\|y^2 U\| \lesssim \|HU\| + \|\Delta U\| + \|U\|$. Therefore,

$$\begin{aligned}\|y_j W_j \tilde{U}\| &\lesssim M_1 e^{(\beta-1/2)\tau} \|\Delta \tilde{U}\|^{n/2q} \big(\|H\tilde{U}\|^{1/2} + \|\Delta \tilde{U}\|^{1/2} + 1\big) \\ &\lesssim M_1 e^{(\beta-1/2)\tau} \big(\|H\tilde{U}\|^{n/2q} + 1\big)\big(\|H\tilde{U}\|^{1/2} + 1\big) \\ &\lesssim M_1 e^{(\beta-1/2)\tau} (\|H\tilde{U}\|^\beta + 1),\end{aligned}$$

where we used $\beta = n/2q + 1/2$. By (3.27) we then obtain

$$\|y_j W_j \tilde{U}\| \lesssim M_1 e^{(\beta-1/2)\tau} \Big(e^{-\beta\tau/2}|a|^\beta (\bar{Q}(\tau))_+^{\beta/2} + e^{-\beta\tau/2}|a|^\beta + 1\big) + \|(A(\tau)-\bar{Q}(\tau)I)\tilde{U}\|^\beta + \bar{Q}(\tau)^\beta + 1\Big). \quad (3.32)$$

By (3.28), (3.30), and (3.32), we get an estimate for $I_4$ which reads

$$\begin{aligned}I_4 &\lesssim \|(A(\tau)-\bar{Q}(\tau)I)\tilde{U}\| \Big(e^{(\beta-1)\tau} M_1 \big(e^{-\beta\tau/2}|a|^\beta (\bar{Q}(\tau))_+^{\beta/2} + e^{-\beta\tau/2}|a|^\beta + 1\big) \\ &\qquad + \|(A(\tau)-\bar{Q}(\tau)I)\tilde{U}\|^\beta + \bar{Q}(\tau)^\beta + 1\big)\Big) \\ &\quad + \|(A(\tau)-\bar{Q}(\tau)I)\tilde{U}\| \Big(e^{(\alpha-1)\tau} M_0 \big(e^{-\alpha\tau/2}|a|^\alpha (\bar{Q}(\tau))_+^{\alpha/2} + e^{-\alpha\tau/2}|a|^\alpha + 1\big) \\ &\qquad + \|(A(\tau)-\bar{Q}(\tau)I)\tilde{U}\|^\alpha + \bar{Q}(\tau)^\alpha + 1\big)\Big).\end{aligned}$$

Combining the bounds for $I_1$, $I_2$, $I_3$, and $I_4$, we get

$$\begin{aligned}&\frac{1}{2}\bar{Q}'(\tau) + \|(A(\tau)-\bar{Q}(\tau)I)\tilde{U}\|^2 \\ &\lesssim e^{-\tau/2}|a|\big(\bar{Q}(\tau)_+^{1/2} + e^{-\tau/2}|a| + 1\big) \\ &\quad + e^{(\beta-1)\tau} M_1 \|(A(\tau)-\bar{Q}(\tau)I)\tilde{U}\| \Big(e^{-\beta\tau/2}|a|^\beta \bar{Q}(\tau)_+^{\beta/2} + e^{-\beta\tau}|a|^{2\beta} + e^{-\beta\tau/2}|a|^\beta \\ &\qquad\qquad + \|(A(\tau)-\bar{Q}(\tau)I)\tilde{U}\|^\beta + \bar{Q}(\tau)^\beta + 1\Big) \\ &\quad + e^{(\alpha-1)\tau} M_0 \|(A(\tau)-\bar{Q}(\tau)I)\tilde{U}\| \Big(e^{-\alpha\tau/2}|a|^\alpha \bar{Q}(\tau)_+^{\alpha/2} + e^{-\alpha\tau}|a|^{2\alpha} + e^{-\alpha\tau/2}|a|^\alpha \\ &\qquad\qquad + \|(A(\tau)-\bar{Q}(\tau)I)\tilde{U}\|^\alpha + \bar{Q}(\tau)^\alpha + 1\Big).\end{aligned} \quad (3.33)$$

We now apply Young's inequality to terms involving $\|(A(\tau)-\bar{Q}(\tau)I)\tilde{U}\|$ on the right-hand side so that we can absorb them into the second term in the left hand side of (3.33). Namely, we use

$$N\|(A(\tau)-\bar{Q}(\tau)I)\tilde{U}\|^\gamma \leq \epsilon_0 \|(A(\tau)-\bar{Q}(\tau)I)\tilde{U}\|^2 + C_{\epsilon_0} N^{2/(2-\gamma)},$$

where $\epsilon_0 \in (0,1]$ is arbitrarily small, with $\gamma = 1, \alpha+1, \beta+1$ and corresponding expressions for $N$. Thus, we obtain

$$\begin{aligned}\bar{Q}'(\tau) &+ \|(A(\tau)-\bar{Q}(\tau)I)\tilde{U}\|^2 \\ &\lesssim e^{-\tau/2}|a|\bar{Q}(\tau)_+^{1/2} + e^{(\beta-2)\tau} M_1^2 |a|^{2\beta} \bar{Q}(\tau)_+^\beta + e^{(2\beta-2)\tau} M_1^2 \bar{Q}(\tau)_+^{2\beta} \\ &\quad + M_0^2 e^{(\alpha-2)\tau}|a|^{2\alpha} \bar{Q}(\tau)_+^\alpha + M_0^2 e^{(2\alpha-2)\tau} \bar{Q}(\tau)_+^{2\alpha} + e^{-\tau}|a|^2 + e^{-\tau/2}|a| \\ &\quad + M_1^2 e^{-2\tau}|a|^{4\beta} + M_1^2 e^{(\beta-2)\tau}|a|^{2\beta} + M_1^{2/(1-\beta)} e^{-2\tau} + M_1^2 e^{(2\beta-2)\tau} \\ &\quad + M_0^2 e^{-2\tau}|a|^{4\alpha} + M_0^2 e^{(\alpha-2)\tau}|a|^{2\alpha} + M_0^{2/(1-\alpha)} e^{-2\tau} + M_0^2 e^{(2\alpha-2)\tau}.\end{aligned} \quad (3.34)$$

Since $|a| \lesssim 1/\epsilon$, the last inequality implies

$$\begin{aligned}\bar{Q}'(\tau) &\lesssim e^{-\tau/2}\epsilon^{-1}\bar{Q}(\tau)_+^{1/2} + e^{(\beta-2)\tau}M_1^2\epsilon^{-2\beta}\bar{Q}(\tau)_+^\beta + e^{(2\beta-2)\tau}M_1^2\bar{Q}(\tau)_+^{2\beta}\\&+ M_0^2 e^{(\alpha-2)\tau}\epsilon^{-2\alpha}\bar{Q}(\tau)_+^\alpha + M_0^2 e^{(2\alpha-2)\tau}\bar{Q}(\tau)_+^{2\alpha} + e^{-\tau}\epsilon^{-2} + e^{-\tau/2}\epsilon^{-1}\\&+ M_1^2 e^{-2\tau}\epsilon^{-4\beta} + M_1^2 e^{(\beta-2)\tau}\epsilon^{-2\beta} + M_1^{2/(1-\beta)}e^{-2\tau} + M_1^2 e^{(2\beta-2)\tau}\\&+ M_0^2 e^{-2\tau}\epsilon^{-4\alpha} + M_0^2 e^{(\alpha-2)\tau}\epsilon^{-2\alpha} + M_0^{2/(1-\alpha)}e^{-2\tau} + M_0^2 e^{(2\alpha-2)\tau}.\end{aligned} \qquad (3.35)$$

For the initial value, we have

$$\begin{aligned}\bar{Q}(\tau_0) &\lesssim Q(\tau_0) + \left|e^{-\tau_0/2}a_j \int y_j \tilde{U}^2\, dy\right| \lesssim 2\epsilon q_{\mathrm{D}}(-\epsilon) + \epsilon^{-1/2}\frac{(HU(\tau_0), U(\tau_0)^{1/2}}{(U(\tau_0), U(\tau_0))^{1/2}}\\&\lesssim 2\epsilon q_{\mathrm{D}}(-\epsilon) + \epsilon^{-1/2}\left(\epsilon q_{\mathrm{D}}(-\epsilon)\right)^{1/2} \lesssim \epsilon q_0 + q_0^{1/2},\end{aligned} \qquad (3.36)$$

which implies

$$\bar{Q}(\tau_0) \leq C_1,$$

where $C_1 = C_0(q_0 + q_0^{1/2})$ and $C_0 \geq 1$ is a constant from (3.36), which is further increased below. Denote by $C_2 \geq 1$ the implicit constant in the inequality (3.35) and choose

$$\epsilon = \frac{1}{2(M_0^a + M_1^b)}, \qquad (3.37)$$

where

$$a = \frac{2}{3 - 4\alpha} \qquad (3.38)$$

and

$$b = \frac{2}{3 - 4\beta}. \qquad (3.39)$$

Note that $\epsilon \in (0, 1/2]$ since $M_0 \geq 1$. We claim that

$$\bar{Q}(\tau) < \frac{C_1 + C_2}{\epsilon}, \qquad \tau \geq \tau_0 = -\log\epsilon \qquad (3.40)$$

if we choose $C_0 \geq 1$ sufficiently large. Assume, contrary to the assertion, that there exists $\tau_1$ such that $\bar{Q}(\tau_1) = (C_1 + C_2)/\epsilon$, and assume that $\tau_1$ is the first such time. Integrating (3.35) between $\tau_0$ and $\tau_1$, we arrive at

$$\begin{aligned}\bar{Q}(\tau_1) &\leq \bar{Q}(\tau_0) + 2C_2(C_1+C_2)^{1/2}\epsilon^{-1} + M_1^2 C_2(C_1+C_2)^\beta \epsilon^{2-4\beta} + 2M_1^2 C_2(C_1+C_2)^{2\beta}\epsilon^{2-4\beta}\\&+ M_0^2 C_2(C_1+C_2)^\alpha \epsilon^{2-4\alpha} + 2M_0^2 C_2(C_1+C_2)^{2\alpha}\epsilon^{2-4\alpha} + C_2\epsilon^{-1} + 2C_2\epsilon^{-1/2}\\&+ M_1^2 C_2 \epsilon^{2-4\beta} + M_1^2 C_2 \epsilon^{2-3\beta} + M_1^{2/(1-\beta)}C_2 \epsilon^2 + 2M_1^2 C_2 \epsilon^{2-2\beta} + M_0^2 C_2 \epsilon^{2-4\alpha}\\&+ M_0^2 C_2 \epsilon^{2-3\alpha} + M_0^{2/(1-\alpha)}C_2 \epsilon^2 + 2M_0^2 C_2 \epsilon^{2-2\alpha},\end{aligned} \qquad (3.41)$$

where we also used $e^{-\tau_0} = \epsilon$ and $0 \leq \alpha, \beta \leq 3/4$. The second term on the right-hand side satisfies

$$2C_2(C_1+C_2)^{1/2}\epsilon^{-1} \leq \frac{C_1 + C_2}{100\epsilon},$$

given $C_1$ is sufficiently large compared to $C_2$, which is achieved by increasing $C_0$. For the third term, we have

$$M_1^2 C_2(C_1+C_2)^\beta \epsilon^{2-4\beta} = M_1^2 C_2(C_1+C_2)^\beta \epsilon^{3-4\beta}\epsilon^{-1} \leq C_2(C_1+C_2)^\beta \epsilon^{-1} \leq \frac{C_1+C_2}{100\epsilon},$$

where in the second inequality we used $M_1^2 \epsilon^{3-4\beta} \leq 1$, and this results from (3.37). In the last step, we used that $C_1$ is sufficiently large compared to $C_2$. We proceed similarly for the rest of the terms in (3.41) and obtain

$$\bar{Q}(\tau_1) \leq \bar{Q}(\tau_0) + \frac{C_1+C_2}{4\epsilon} \leq C_1 + \frac{C_1+C_2}{4\epsilon} < \frac{C_1+C_2}{\epsilon},$$

where in the last step we used $\epsilon \leq 1/2$. This is a contradiction with a choice of $\tau_1$, and thus we conclude that (3.40) holds for all $\tau \geq \tau_0$. Finally, by (3.37) and (3.40), we get $\bar{Q}(\tau) \lesssim M_0^a + M_1^b$, as desired. $\square$

Now we show that the modified frequency function $\bar{Q}$ controls the vanishing order of $u$. We first prove the following lemma, which shows the convergence of this function as $\tau \to \infty$ and the order of vanishing of $\int u^2(x,t)G(x,t)\, dx$, where $G$ is defined in (3.1).

LEMMA 3.3. *Under the assumptions of Theorem 2.1, the modified frequency function satisfies $\bar{Q}(\tau) \to m/2$ as $\tau \to \infty$ for some $m \in \mathbb{N}$ such that $m \lesssim M_0^a + M_1^b$, where $a$ and $b$ are as in (3.38) and (3.39). Moreover, for all $\delta > 0$, there exist $\eta \in (\log(1/\epsilon), 0)$ and $A_1(\delta), A_2(\delta) > 0$ such that*

$$A_1(\delta)|t|^{m+\delta} \leq \int u^2(x,t) G(x,t) \, dx \leq A_2(\delta)|t|^{m-\delta}, \tag{3.42}$$

*for all $t \in [-\eta, 0)$.*

PROOF OF LEMMA 3.3. Taking the inner product of (3.15) with $U$, we obtain

$$\frac{1}{2}\frac{d}{d\tau}\|U\|^2 + (A(\tau)U, U) = e^{-\tau/2}(y_j W_j U, U) + e^{-\tau/2}(W_j \partial_j U, U) + e^{-\tau}(VU, U), \tag{3.43}$$

where we used $(e^{-\tau/2} a_j \partial_j U, U) = 0$. Denoting the right-hand side of (3.43) by $f(\tau)$, we obtain

$$\begin{aligned}
f(\tau) &= e^{-\tau/2} \int y_j W_j U^2 \, dy + e^{-\tau/2} \int W_j U \partial_j U \, dy + e^{-\tau} \int V U^2 \, dy \\
&\lesssim e^{(\beta-1)\tau} M_1 \|yU^2\|_{L^{q/(q-1)}(\mathbb{R}^n)} + e^{(\beta-1)\tau} M_1 \|U\nabla U\|_{L^{q/(q-1)}(\mathbb{R}^n)} + e^{(\alpha-1)\tau} M_0 \|U^2\|_{L^{p/(p-1)}(\mathbb{R}^n)} \\
&\lesssim e^{(\beta-1)\tau} M_1 \|yU\| \|U\|_{L^{2q/(q-2)}(\mathbb{R}^n)} + e^{(\beta-1)\tau} M_1 \|U\|_{L^{2q/(q-2)}(\mathbb{R}^n)} \|\nabla U\| \\
&\quad + e^{(\alpha-1)\tau} M_0 \|U\|^2_{L^{2p/(p-1)}(\mathbb{R}^n)}.
\end{aligned}$$

We now use the Gagliardo-Nirenberg inequality to get

$$\begin{aligned}
f(\tau) &\lesssim e^{(\beta-1)\tau} M_1 \|yU\| \|U\|^{1-n/q} \|\nabla U\|^{n/q} + e^{(\beta-1)\tau} M_1 \|U\|^{1-n/q} \|\nabla U\|^{n/q+1} \\
&\quad + e^{(\alpha-1)\tau} M_0 \|U\|^{2-n/p} \|\nabla U\|^{n/p}.
\end{aligned}$$

Therefore, by (3.26), we may estimate

$$\frac{f(\tau)}{\|U\|^2} \lesssim e^{(\beta-1)\tau} M_1 \left(\bar{Q}(\tau)_+ + e^{-\tau}|a|^2 + 1\right)^\beta + e^{(\alpha-1)\tau} M_0 \left(\bar{Q}(\tau)_+ + e^{-\tau}|a|^2 + 1\right)^\alpha, \tag{3.44}$$

where $\alpha$ and $\beta$ are as in (3.23) and (3.29). For convenience, we denote

$$M = M_0^a + M_1^b. \tag{3.45}$$

Now, we integrate (3.44) from $\tau_0$ to $\tau$ and obtain

$$\begin{aligned}
\int_{\tau_0}^{\tau} \frac{f(s)}{\|U(\cdot, s)\|^2} \, ds &\lesssim M^\beta M_1 (e^{-(1-\beta)\tau_0} - e^{-(1-\beta)\tau}) + |a|^{2\beta} M_1 (e^{-\tau_0} - e^{-\tau}) + M_1 (e^{-(1-\beta)\tau_0} - e^{-(1-\beta)\tau}) \\
&\quad + M^\alpha M_0 (e^{-(1-\alpha)\tau_0} - e^{-(1-\alpha)\tau}) + |a|^{2\alpha} M_0 (e^{-\tau_0} - e^{-\tau}) + M_0 (e^{-(1-\alpha)\tau_0} - e^{-(1-\alpha)\tau}),
\end{aligned}$$

where we used (3.16). On the other hand, if we integrate from $\tau_1$ to $\tau$ the equation

$$\frac{1}{2\|U\|^2} \frac{d}{d\tau} \|U\|^2 + \bar{Q}(\tau) = \frac{f(\tau)}{\|U\|^2},$$

we get

$$\begin{aligned}
\frac{1}{2} \log \|U(\cdot, \tau)\|^2 - \frac{1}{2} \log \|U(\cdot, \tau_1)\|^2 &= -\int_{\tau_1}^{\tau} \bar{Q}(s) \, ds + \int_{\tau_1}^{\tau} \frac{f(s)}{\|U(\cdot, s)\|^2} \, ds \\
&\leq -\frac{m}{2}(\tau - \tau_1) - \int_{\tau_1}^{\tau} \left(\bar{Q}(s) - \frac{m}{2}\right) ds + C(e^{-\tau_1/4} + e^{-\tau/4}),
\end{aligned} \tag{3.46}$$

where the constant depends on $|a|$, $M$ (note that $\alpha$ and $\beta$ are considered fixed). To show that $\bar{Q}(\tau) \to m/2$ as $\tau \to \infty$ for some $m \lesssim M$, we use $\bar{Q}(\tau) \lesssim M$ in (3.34) and integrate between $\tau_1$ and $\tau$, where $0 \leq \tau_1 \leq \tau$ obtaining

$$\bar{Q}(\tau) - \bar{Q}(\tau_1) + \int_{\tau_1}^{\tau} \|(A(s) - \bar{Q}(s)I)\tilde{U}\|^2 \, ds \lesssim e^{-\tau_1/4} + e^{-\tau/4}.$$

Therefore,

$$\int_{\tau_0}^{\tau} \|(A(s) - \bar{Q}(s)I)\tilde{U}\|^2 \, ds \lesssim M + e^{-\tau_1/4} + e^{-\tau/4}.$$

We also have

$$\|(A(\tau) - \bar{Q}(\tau)I)\tilde{U}\|^2 \gtrsim \frac{1}{2}\|(H(\tau) - \bar{Q}(\tau)I)\tilde{U}\|^2 - e^{-\tau}|a|^2 \bar{Q} - e^{-\tau}|a|^2$$

by definitions of $A(\tau)$ and $H(\tau)$. Therefore,

$$\int_{\tau_0}^{\tau} \|(H(s) - \bar{Q}(s)I)\tilde{U}\|^2 \, ds \lesssim \int_{\tau_0}^{\tau} \|(A(s) - \bar{Q}(s)I)\tilde{U}\|^2 \, ds + e^{-\tau} \lesssim M + e^{-\tau_0/2} + e^{-\tau/2} < \infty, \qquad (3.47)$$

where we used the boundedness of $\bar{Q}$ and $a$ in the first inequality. Combining (3.47) with $\|(H - \bar{Q}(\tau)I)\tilde{U}\| \geq \operatorname{dist}(\bar{Q}(\tau), \operatorname{sp}(H))$ and recalling that $\|\tilde{U}\| = 1$, we conclude that $\operatorname{dist}(\bar{Q}(\tau), \operatorname{sp}(H)) \to 0$ as $\tau \to \infty$. It is well-known that $\operatorname{sp}(H) = \{m/2 : m \in \mathbb{N}_0\}$, (cf. [CK, p.664]). Thus we obtain $\bar{Q}(\tau) \to m/2$ as $\tau \to \infty$ for some $m \lesssim M_0^a + M_1^b$.

Now, from (3.46), we deduce that for all $\delta > 0$ there exists $\tau_1 > 0$ depending on $\delta$ such that

$$-\delta(\tau - \tau_1) \leq \log \|U(\tau)\|^2 - \log \|U(\tau_1)\|^2 + m(\tau - \tau_1) \leq \delta(\tau - \tau_1), \qquad \tau \geq \tau_1.$$

Then,

$$e^{-\delta(\tau - \tau_1)} \leq \frac{e^{\tau m}\|U(\tau)\|^2}{e^{\tau_1 m}\|U(\tau_1)\|^2} \leq e^{\delta(\tau - \tau_1)}.$$

Therefore, there exists $A_1(\delta), A_2(\delta) > 0$ such that

$$A_1(\delta)e^{-\delta\tau} \leq e^{\tau m}\|U(\tau)\|^2 \leq A_2(\delta)e^{\delta\tau}. \qquad (3.48)$$

Restating (3.48) in the original variable $u$, we obtain (3.42). $\square$

The following lemma provides control on $\int_{B(0,R)} Q(x,t)G(x,t) \, dx$, with $Q$ a homogeneous polynomial of degree $d$. We use this lemma in our proof of Theorem 2.1.

LEMMA 3.4. *Let $Q = \Sigma_{|\mu|+2l=d} C_{\mu,l} x^{\mu} t^l$ be a homogeneous polynomial of degree $d \in \mathbb{N}$. Then,*

$$\int_{\mathbb{R}^n} Q(x,t)G(x,t) \, dx \lesssim |t|^{d/2}, \qquad (3.49)$$

*where the constant in (3.49) depends on the polynomial only. Moreover, if all the coordinates of $\mu = (\mu_1, \mu_2, \ldots, \mu_n)$ are even, then for all $R > 0$,*

$$\int_{B(0,R)} x^{\mu} t^l G(x,t) \, dx \lesssim |t|^{l+|\mu|/2},$$

*as $t \to 0^-$. If $\mu_i$ is an odd integer for some $i \in \{1, \ldots, n\}$, then*

$$\int_{B(0,R)} x^{\mu} t^l G(x,t) \, dx = 0.$$

For the proof of Lemma 3.4, cf. [CK, p.670]).

PROOF OF THEOREM 2.1. Without loss of generality, let $I = [0, T]$ and $(x_0, t_0) = (0, 0)$. Let $m$ be as in the statement of Lemma 3.3, and denote by $d$ the vanishing order of $u$ at $(0, 0)$. We claim that $d \leq m$. Since the degree of vanishing of $u$ at $(0,0)$ is $d$, we have

$$|u(x,t)| \lesssim (|x|^2 + |t|)^{d/2},$$

for all $(x,t) \in Q_1(0,0)$ with $Q_1(0,0)$ defined in (2.4). Also, for any $R > 0$, we have

$$\int_{\mathbb{R}^n \setminus B(0,R)} u^2(x,t) G(x,t) \, dx \leq \|u\|_{L^\infty(\Omega)}^2 \int_R^\infty \int_{\partial B(0,R)} \frac{e^{-r^2/4|t|}}{|t|^{n/2}} \, dS_y \, dr \lesssim \|u\|_{L^\infty(\Omega)}^2 \frac{|t|^{1/2}}{R} e^{-R^2/8|t|}.$$

By Lemma 3.3, there exists $\eta \in (0, \epsilon)$ and $A_1(\delta), A_2(\delta) > 0$ such that

$$A_1(\delta)|t|^{m+\delta} \leq u^2(x,t) G(x,t) \, dx \leq A_2(\delta)|t|^{m-\delta}, \qquad t \in [-\eta, 0]. \qquad (3.50)$$

Moreover, by Lemma 3.4,

$$\int_{B(0,R)} u^2(x,t) G(x,t) \, dx \lesssim \int_{B(0,R)} (|x|^2 + |t|)^d G(x,t) \, dx \lesssim |t|^d, \qquad (3.51)$$

where the constant in (3.51) depends on $R$. Combining (3.50) with (3.51), we get, for a fixed $\delta \in (0, 1]$,

$$A_1(\delta)|t|^{m+\delta} \lesssim |t|^d, \qquad t \in [-\eta, 0].$$

Sending $t \to 0$, we obtain

$$d \leq m + \delta,$$

and then letting $\delta \to 0$, we conclude that $d \leq m \lesssim M_0^a + M_0^b$ as desired.

Therefore, we have proven Theorem 2.1 when $v$, $w$, and $u_0$ smooth. For the general case, we approximate $v$, $w$, and $u_0$ by smooth functions $v^\epsilon$, $w^\epsilon$, and $u_0^\epsilon$ in $L_t^\infty L_x^p$, $L_t^\infty L_x^q$, and $L_t^\infty L_x^\infty$, respectively. Let $u^\epsilon$ be a solution of the equation

$$\begin{aligned}\partial_t u^\epsilon - \Delta u^\epsilon &= w^\epsilon \partial_j u^\epsilon + v^\epsilon u^\epsilon \\ u^\epsilon(x, 0) &= u_0^\epsilon(x).\end{aligned} \tag{3.52}$$

Subtracting (3.52) from (2.1), we get

$$\begin{aligned}\partial_t \tilde{u} - \Delta \tilde{u} &= w^\epsilon \partial_j \tilde{u} + \tilde{w} \partial_j u + v^\epsilon \tilde{u} + \tilde{v} u \\ \tilde{u}(x, 0) &= u_0^\epsilon(x) - u_0(x),\end{aligned} \tag{3.53}$$

where $\tilde{u} = u^\epsilon - u$, $\tilde{v} = v^\epsilon - v$ and $\tilde{w} = w^\epsilon - w$. First, we have $u, \tilde{u} \in L^\infty L^2 \cap L^2 H^1$. Using $v \in L^\infty L^{\max\{2,n/2\}}$ and $w \in L^\infty L^{2n}$, bootstrapping, and the $W^{2,p}$ maximal regularity for parabolic equations, we get

$$D^2 u, D^2 \tilde{u} \in L^2 L^{2n}(\Omega \times (0, T)). \tag{3.54}$$

Taking the inner product of (3.53) with $\tilde{u}$, we obtain

$$\begin{aligned}\frac{1}{2}\frac{d}{dt}\|\tilde{u}(\cdot, t)\|_{L^2(\Omega)}^2 + \|\nabla \tilde{u}\|_{L^2(\Omega)}^2 &= \int w^\epsilon \tilde{u} \partial_j \tilde{u} + \int \tilde{w} \tilde{u} \partial_j u + \int v^\epsilon \tilde{u}^2 + \int \tilde{v} u \tilde{u} \\ &\lesssim \|w^\epsilon\|_{L^q(\Omega)}\|\tilde{u}\|_{L^{2q/(q-2)}(\Omega)}\|\nabla \tilde{u}\|_{L^2(\Omega)} + \|\tilde{w}\|_{L^q(\Omega)}\|\tilde{u}\|_{L^{2q/(q-2)}(\Omega)}\|\nabla u\|_{L^2(\Omega)} \\ &\quad + \|v^\epsilon\|_{L^p(\Omega)}\|\tilde{u}\|_{L^{2p/(p-1)}(\Omega)}^2 + \|\tilde{v}\|_{L^p(\Omega)}\|u\|_{L^{2p/(p-1)}(\Omega)}\|\tilde{u}\|_{L^{2p/(p-1)}(\Omega)}.\end{aligned}$$

Using $\|\tilde{u}\|_{L^{2q/(q-2)}(\Omega)} \lesssim \|\tilde{u}\|_{L^2(\Omega)}^{1-n/q}\|\nabla \tilde{u}\|_{L^2(\Omega)}^{n/q}$ and $\|\tilde{u}\|_{L^{2p/(p-1)}(\Omega)} \lesssim \|\tilde{u}\|_{L^2(\Omega)}^{1-n/2p}\|\nabla \tilde{u}\|_{L^2(\Omega)}^{n/2p}$, with the same inequality for $u$, and applying Young's inequality, we get, for an arbitrary $\epsilon_0 \in (0, 1]$,

$$\begin{aligned}\frac{d}{dt}\|\tilde{u}(\cdot, t)\|_{L^2(\Omega)}^2 &\leq C\|w^\epsilon\|_{L^q(\Omega)}^{2q/(q-n)}\|\tilde{u}\|_{L^2(\Omega)}^2 + C_{\epsilon_0}\|\tilde{w}\|_{L^q(\Omega)}^{2q/(q-n)}\|\tilde{u}\|_{L^2(\Omega)}^2 \\ &\quad + C\|v^\epsilon\|_{L^p(\Omega)}^{2p/(2p-n)}\|\tilde{u}\|_{L^2(\Omega)}^2 + C_{\epsilon_0}\|\tilde{v}\|_{L^p(\Omega)}^{2p/(2p-n)}\|\tilde{u}\|_{L^2(\Omega)}\|u\|_{L^2(\Omega)} + \epsilon_0\|\nabla u\|_{L^2(\Omega)}^2.\end{aligned} \tag{3.55}$$

Integrating (3.55) from 0 to $t_0$, we obtain

$$\begin{aligned}&\|\tilde{u}(\cdot, t_0)\|_{L^2(\Omega)}^2 - \|\tilde{u}(\cdot, 0)\|_{L^2(\Omega)}^2 \\ &\leq C\int_0^{t_0} \|w^\epsilon\|_{L^q(\Omega)}^{2q/(q-n)(\Omega)}\|\tilde{u}\|_{L^2(\Omega)}^2 + C_{\epsilon_0}\int_0^{t_0} \|\tilde{w}\|_{L^q(\Omega)}^{2q/(q-n)}\|\tilde{u}\|_{L^2(\Omega)}^2 \\ &\quad + C\int_0^{t_0} \|v^\epsilon\|_{L^p(\Omega)}^{2p/(2p-n)}\|\tilde{u}\|_{L^2(\Omega)}^2 + C_{\epsilon_0}\int_0^{t_0} \|\tilde{v}\|_{L^p(\Omega)}^{2p/(2p-n)}\|\tilde{u}\|_{L^2(\Omega)}\|u\|_{L^2(\Omega)} + \epsilon_0\int_0^{t_0}\|\nabla u\|_{L^2(\Omega)}^2,\end{aligned}$$

for every $\epsilon_0 \in (0, 1]$. Sending $\epsilon \to 0^+$, it follows that

$$\limsup_{\epsilon \to 0} \|\tilde{u}(\cdot, t_0)\|_{L^2(\Omega)}^2 \leq \epsilon_0 \int_0^{t_0} \|\nabla u\|_{L^2(\Omega)}^2, \qquad t_0 \geq 0,$$

for every $\epsilon_0 \in (0, 1]$. Sending then $\epsilon_0 \to 0$, we then obtain

$$\|u^\epsilon(\cdot, t_0)\|_{L^2(\Omega)} \to \|u(\cdot, t_0)\|_{L^2(\Omega)}, \qquad t_0 \geq 0. \tag{3.56}$$

Now, we claim that $\nabla u^\epsilon \to \nabla u$ in $L_t^\infty L_x^2$. To prove this, we test the first equation in (3.53) with $-\Delta \tilde{u}$, which we may by (3.54), and write

$$\frac{d}{dt}\|\nabla \tilde{u}(\cdot, t)\|_{L^2(\Omega)}^2 + \|\Delta \tilde{u}(\cdot, t)\|_{L^2(\Omega)}^2 = -\int w^\epsilon \partial_j \tilde{u} \Delta \tilde{u}\, dx - \int \tilde{w} \partial_j u \Delta \tilde{u}\, dx - \int v^\epsilon \tilde{u} \Delta \tilde{u}\, dx - \int \tilde{v} u \Delta \tilde{u}\, dx.$$

By Hölder's and Sobolev's inequalities, we have

$$\begin{aligned}
&\frac{d}{dt}\|\nabla \tilde{u}(\cdot,t)\|^2_{L^2(\Omega)} + \|\Delta \tilde{u}(\cdot,t)\|^2_{L^2(\Omega)}\\
&\lesssim \|w^\epsilon\|_{L^q(\Omega)}\|\nabla\tilde{u}\|_{L^{2q/(q-2)}(\Omega)}\|\Delta\tilde{u}\|_{L^2(\Omega)} + \|\tilde{w}\|_{L^q(\Omega)}\|\nabla u\|_{L^{2q/(q-2)}(\Omega)}\|\Delta\tilde{u}\|_{L^2(\Omega)}\\
&\quad + \|v^\epsilon\|_{L^p(\Omega)}\|\tilde{u}\|_{L^{2p/(p-2)}(\Omega)}\|\Delta\tilde{u}\|_{L^2(\Omega)} + \|\tilde{v}\|_{L^p(\Omega)}\|u\|_{L^{2p/(p-2)}(\Omega)}\|\Delta\tilde{u}\|_{L^2(\Omega)}\\
&\lesssim \|w^\epsilon\|_{L^q(\Omega)}\|\tilde{u}\|_{L^2(\Omega)}^{1/2-n/2q}\|\Delta\tilde{u}\|^{3/2+n/2q}_{L^2(\Omega)} + \|\tilde{w}\|_{L^q(\Omega)}\|u\|^{1/2-n/2q}\|\Delta u\|^{1/2+n/2q}_{L^2}\|\Delta\tilde{u}\|_{L^2(\Omega)}\\
&\quad + \|v^\epsilon\|_{L^p(\Omega)}\|\tilde{u}\|^{1-n/2p}_{L^2(\Omega)}\|\Delta\tilde{u}\|^{n/2p+1}_{L^2(\Omega)} + \|\tilde{v}\|_{L^p(\Omega)}\|u\|^{1-n/2p}_{L^2(\Omega)}\|\Delta u\|^{n/2p}_{L^2(\Omega)}\|\Delta\tilde{u}\|_{L^2(\Omega)}.
\end{aligned} \quad (3.57)$$

Applying Young's inequality on all four terms and absorbing the factors containing $\Delta\tilde{u}$, we get

$$\begin{aligned}
\frac{d}{dt}\|\nabla\tilde{u}(\cdot,t)\|^2_{L^2(\Omega)} &\lesssim \|w^\epsilon\|^{4q/(q-n)}_{L^q(\Omega)}\|\tilde{u}\|^2_{L^2(\Omega)} + \|\tilde{w}\|^2_{L^q(\Omega)}\|u\|^{1-n/q}_{L^2(\Omega)}\|\Delta u\|^{1+n/q}_{L^2(\Omega)}\\
&\quad + \|v^\epsilon\|^{4p/(2p-n)}_{L^p(\Omega)}\|\tilde{u}\|^2_{L^2(\Omega)} + \|\tilde{v}\|^2_{L^p(\Omega)}\|u\|^{2-n/p}_{L^2(\Omega)}\|\Delta u\|^{n/p}_{L^2(\Omega)}.
\end{aligned} \quad (3.58)$$

For any $t_0 \geq 0$, we integrate (3.58) from 0 to $t_0$ to obtain

$$\begin{aligned}
&\|\nabla\tilde{u}(\cdot,t_0)\|^2_{L^2(\Omega)} - \|\nabla\tilde{u}(\cdot,0)\|^2_{L^2(\Omega)}\\
&\lesssim t_0\|w^\epsilon\|_{L^\infty L^q(\Omega\times I)}\|\tilde{u}\|^2_{L^\infty L^2(\Omega\times I)} + \|\tilde{w}\|^2_{L^\infty L^q(\Omega\times I)}\|u\|^{1-n/q}_{L^\infty L^2(\Omega\times I)}\int_0^{t_0}\|\Delta\tilde{u}\|^{1+n/q}_{L^2(\Omega)}\\
&\quad + t_0\|v^\epsilon\|^{4p/(2p-n)}_{L^\infty L^p(\Omega\times I)}\|\tilde{u}\|^2_{L^\infty L^2(\Omega\times I)} + \|\tilde{v}\|^2_{L^\infty L^p(\Omega\times I)}\|u\|^{2-n/p}_{L^2 L^{2p/(p-2)}(\Omega\times I)}\int_0^{t_0}\|\Delta u\|^{n/p}_{L^2 L^2(\Omega\times I)},
\end{aligned}$$

where $\|\Delta\tilde{u}\|_{L^2 L^2(\Omega\times I)} < \infty$ due to (3.54). Letting $\epsilon \to 0^+$, we obtain

$$\|\nabla u^\epsilon(\cdot,t_0)\|_{L^2(\Omega)} \to \|\nabla u(\cdot,t_0)\|_{L^2(\Omega)}, \qquad t_0 \geq 0. \quad (3.59)$$

Denote by $\bar{Q}^\epsilon(\tau)$ and $\bar{Q}(\tau)$ the modified frequency functions corresponding to $u^\epsilon$ and $u$, respectively. Recalling that

$$\bar{Q}(\tau) = \frac{(HU,U)}{\|U\|^2} - \frac{e^{-\tau/2}a_j}{\|U\|^2}\int y_j U^2\, dy,$$

we now claim that $\bar{Q}^\epsilon(\tau) \to \bar{Q}(\tau)$ as $\epsilon \to 0$ for all $\tau$. Note that

$$\|U^\epsilon(\tau)\|^2 - \|U(\tau)\|^2 = \int(u_\epsilon^2 - u^2)G(x,t)\,dx \lesssim \|\tilde{u}\|_{L^2(\Omega)}\|(u^\epsilon + u)\|_{L^2(\Omega)} \to 0, \quad (3.60)$$

as $\epsilon \to 0^+$ by (3.56). For $(HU_\epsilon, U_\epsilon) - (HU,U)$, we have

$$(HU_\epsilon, U_\epsilon) - (HU,U) = |t|\int(|\nabla u_\epsilon|^2 - |\nabla u|^2)G(x,t)\,dx \lesssim |t|\|\nabla\tilde{u}\|_{L^2(\Omega)}\|\nabla(u^\epsilon + u)\|_{L^2(\Omega)} \to 0,$$

as $\epsilon \to 0^+$ by (3.59). Lastly, we estimate $\int y_j U_\epsilon^2\,dy - \int y_j U^2\,dy$ as

$$\begin{aligned}
\int y_j U_\epsilon^2\,dy - \int y_j U^2\,dy &= \int y_j(U_\epsilon^2 - U^2) \lesssim \|y(U_\epsilon - U)\|\|U_\epsilon + U\|\\
&\lesssim \big((H(U_\epsilon - U), U_\epsilon - U)^{1/2} + \|U_\epsilon - U\|\big)\|U^\epsilon + U\|\\
&\lesssim \big((|t|\|\nabla(u^\epsilon - u)\|\|\nabla(u^\epsilon - u)G\|)^{1/2} + \|U_\epsilon - U\|\big)\|U_\epsilon + U\| \to 0,
\end{aligned} \quad (3.61)$$

as $\epsilon \to 0^+$ where we used $\|yU\|^2 \lesssim (HU,U) + \|U\|^2$ in the second inequality. By (3.60)–(3.61), we obtain $\bar{Q}^\epsilon(\tau) \to \bar{Q}(\tau)$ as $\epsilon \to 0$ for all $\tau$. Note that we have

$$\bar{Q}^\epsilon(\tau) \lesssim M_0^a + M_1^b, \quad (3.62)$$

uniformly in $\tau$ and $\epsilon \in (0,1]$. Passing to the limit, we get

$$\bar{Q}(\tau) \lesssim M_0^a + M_1^b, \qquad \tau \geq 0. \quad (3.63)$$

Next, using $\operatorname{dist}(\bar{Q}^\epsilon(\tau),\operatorname{sp}(H)) \lesssim \|(H - \bar{Q}(\tau)I)\tilde{U}\|$ and (3.47), we obtain, for all $\epsilon \in (0,1]$

$$\int_{\tau_0}^\tau \operatorname{dist}(\bar{Q}^\epsilon(\tau),\operatorname{sp}(H))^2 \lesssim M_0^a + M_1^b + e^{-\tau_0/2} + e^{-\tau/2} \lesssim M_0^a + M_1^b = K, \quad (3.64)$$

where we used $M_0, M_1 \geq 1$. Observe that by (3.62) we have $Q^\epsilon(\tau) \in [0, \bar{M}]$, where $\bar{M}$ is independent of $\epsilon$ and $\tau$. Thus letting $\epsilon \to 0$ in (3.64) yields

$$\int_{\tau_0}^{\tau} \text{dist}(\bar{Q}(\tau), \text{sp}(H))^2 \leq K. \tag{3.65}$$

Combining (3.63), (3.65), and $\text{sp}(H) = \{m/2 : m \in \mathbb{N}_0\}$, we get $\bar{Q}(\tau) \to m/2$ for some $m \lesssim M_0^a + M_1^b$ with $a$ and $b$ as in (2.7). The rest of the proof is similar to the case when $v$, $w$, and $u_0$ are smooth. $\square$

## 4. Doubling property

In the proof of Theorem 2.2, we need the following observability lemmas.

LEMMA 4.1. *Let $K \geq 0$ and $\delta_0 \in (0, 1/2]$. If a $\Omega$-periodic function $u$ satisfies*

$$\|u\|_{L^2(\Omega)}^2 \leq e^K \int_{\mathbb{R}^n} u(x)^2 G(x, -\delta^2) \, dx$$

*with $\delta \in (0, \delta_0]$ so small that*

$$\frac{1}{\delta^2} \geq \frac{C}{\delta_0^2} \log \frac{1}{\delta} + \frac{C(K+1)}{\delta_0^2},$$

*for a sufficiently large constant $C > 0$, then*

$$\|u\|_{L^2(\Omega)}^2 \leq \frac{Ce^K}{\delta^n} \|u\|_{L^2(B_{\delta_0})}^2,$$

*where $\Omega = [0, 1]^n$.*

For the proof, cf. [Ku4, p.775].

LEMMA 4.2. *Under assumptions of Theorem 2.1, let $t_0 < 0$ such that $|t_0| < \epsilon$ with $\epsilon$ defined in (3.37). Then, for any $\delta_0 \in (0, 1/2]$, we have*

$$\|u(\cdot, t)\|_{L^2(\Omega)}^2 \lesssim e^{(n+1)\delta_0^2/\delta^2} \|u(\cdot, t)\|_{L^2(B_{\delta_0})}^2, \qquad t \in [-\delta^2, 0] \tag{4.1}$$

*with $\delta > 0$ satisfying*

$$\frac{1}{\delta^2} \geq \frac{\log 1/\delta}{\delta_0^2} + \frac{\gamma(\delta) + 1}{\delta_0^2} \tag{4.2}$$

*and*

$$\delta < \min\{\sqrt{\epsilon}, \delta_0\}, \tag{4.3}$$

*where*

$$\gamma(\delta) = M_1^2 + M_0 + M_1 M^{2\beta-1} + M_0 M^{2\alpha-1} + M \log(1/\delta^2), \tag{4.4}$$

*with $M$ defined in (3.45). Moreover,*

$$\delta = \min\left\{\sqrt{\epsilon}, \frac{\delta_0}{\log(1/\delta_0)}, \frac{\delta_0^2}{M^2}, \frac{\delta_0}{M_1 + M_0^{1/2} + M^{\beta-1/2}M_0^{1/2} + M_0^{1/2}M^{\alpha-1/2}}\right\}$$

*satisfies the conditions (4.2) and (4.3). Therefore, $(n+1)\delta_0^2/\delta^2$ can be replaced by a polynomial $P$ with respect $\delta_0$ and $P$ depends on $n$, $M_0$, $M_1$, $p$, and $q$.*

PROOF OF LEMMA 4.2. From (3.44), we see that

$$\frac{f(\tau)}{\|U\|^2} \gtrsim -e^{(\beta-1)\tau} M_1 (\bar{Q}(\tau)_+ + e^{-\tau}|a|^2 + 1)^\beta - e^{(\alpha-1)\tau} M_0 (\bar{Q}(\tau)_+ + e^{-\tau}|a|^2 + 1)^\alpha,$$

where $\alpha$ and $\beta$ are as in (3.23) and (3.29). Therefore, by replacing $\tau_1$ by $\tau_0$ in (3.46), we have

$$\frac{1}{2}\log\|U(\cdot, \tau)\|^2 - \frac{1}{2}\log\|U(\cdot, \tau_0)\|^2 = -\int_{\tau_0}^{\tau} \bar{Q}(s)\,ds + \int_{\tau_0}^{\tau} \frac{f(s)}{\|U(\cdot, s)\|^2}\,ds$$

$$\gtrsim -(\tau - \tau_0)\frac{2C_1 + C_2}{\epsilon} - e^{(\beta-1)\tau_0} M_1\Big(\frac{2C_1 + C_2}{\epsilon} + e^{-\tau_0}|a|^2 + 1\Big)^\beta$$

$$- e^{(\alpha-1)\tau_0} M_0\Big(\frac{2C_1 + C_2}{\epsilon} + e^{-\tau_0}|a|^2 + 1\Big)^\alpha,$$

where we used $\bar{Q}(\tau) \lesssim (2C_1 + C_2)/\epsilon$ for all $\tau \geq \tau_0$ in the last inequality. Since $\epsilon = e^{-\tau_0}$, we have

$$\log \|U(\cdot, \tau)\|^2 - \log \|U(\cdot, \tau_0)\|^2 \gtrsim -(\tau - \tau_0)\epsilon^{-1} - \epsilon^{1-\beta}M_1(\epsilon^{-1} + \epsilon|a|^2 + 1)^\beta - \epsilon^{1-\alpha}M_0(\epsilon^{-1} + \epsilon|a|^2 + 1)^\alpha. \quad (4.5)$$

Let $q_\mathrm{D}(t) = \|\nabla u\|^2_{L^2(\Omega)}/\|u\|^2_{L^2(\Omega)}$. We now derive a similar estimate for $\|u(\cdot, t)\|^2_{L^2(\Omega)}/\|u(\cdot, -\epsilon)\|^2_{L^2(\Omega)}$. Taking the inner product of (3.7) with $u$ and using $\|\nabla u\|^2_{L^2(\Omega)} = q_\mathrm{D}(t)\|u\|^2_{L^2(\Omega)}$, we obtain

$$\begin{aligned}\frac{d}{dt}\|u(\cdot,t)\|^2_{L^2(\Omega)} &+ q_\mathrm{D}(t)\|u(\cdot,t)\|^2_{L^2(\Omega)} \\ &\lesssim \|w(\cdot,t)\|_{L^q(\Omega)}\|u(\cdot,t)\|_{L^{2q/(q-2)}(\Omega)}\|u(\cdot,t)\|_{L^2(\Omega)} q_\mathrm{D}(t)^{1/2} + \|v(\cdot,t)\|_{L^p(\Omega)}\|u(\cdot,t)\|^2_{L^{2p/(p-1)}(\Omega)} \\ &\lesssim M_1 \|u(\cdot,t)\|^{2-2\beta}_{L^2(\Omega)}\|\nabla u(\cdot,t)\|^{2\beta-1}_{L^2(\Omega)}\|u(\cdot,t)\|_{L^2(\Omega)} q_\mathrm{D}(t)^{1/2} + M_0\|u(\cdot,t)\|^{2-2\alpha}_{L^2(\Omega)}\|\nabla u(\cdot,t)\|^{2\alpha}_{L^2(\Omega)}\end{aligned} \quad (4.6)$$

since $\int_\Omega a_j u \partial_j u \, dx$ vanishes. Dividing both sides of (4.6) by $\|u(\cdot,t)\|^2_{L^2(\Omega)}$, we get

$$\begin{aligned}\frac{1}{\|u(\cdot,t)\|^2_{L^2(\Omega)}}\frac{d}{dt}\|u(\cdot,t)\|^2_{L^2(\Omega)} + q_\mathrm{D}(t) &\lesssim M_1 \frac{\|\nabla u(\cdot,t)\|^{2\beta-1}_{L^2(\Omega)}}{\|u(\cdot,t)\|^{2\beta-1}_{L^2(\Omega)}} q_\mathrm{D}(t)^{1/2} + M_0 \frac{\|\nabla u(\cdot,t)\|^{2\alpha}_{L^2(\Omega)}}{\|u(\cdot,t)\|^{2\alpha}_{L^2(\Omega)}} \\ &\lesssim M_1 q_0^{\beta-1/2} q_\mathrm{D}(t)^{1/2} + M_0 q_0^\alpha \lesssim M_1 q_\mathrm{D}(t)^{1/2} + M_0\end{aligned} \quad (4.7)$$

if we allow the last implicit constant to depend on $q_0$. Applying Young's inequality to the right-hand side of (4.7), we obtain

$$\frac{1}{\|u(\cdot,t)\|^2_{L^2(\Omega)}}\frac{d}{dt}\|u(\cdot,t)\|^2_{L^2(\Omega)} \lesssim M_1^2 + M_0. \quad (4.8)$$

Integrating (4.8) from $-\epsilon$ to $t$ leads to

$$\log \frac{\|u(\cdot,t)\|^2_{L^2(\Omega)}}{\|u(\cdot,-\epsilon)\|^2_{L^2(\Omega)}} \lesssim (M_1^2 + M_0)(t + \epsilon) \lesssim (M_1^2 + M_0)\epsilon, \qquad t \in (\epsilon, 0). \quad (4.9)$$

By definition of $U$ in (3.9), we have

$$\frac{\|U(\cdot,\tau_0)\|^2}{\|u(\cdot,-\epsilon)\|^2_{L^2(\Omega)}} = \frac{\int_{\mathbb{R}^n} u(x,-\epsilon)^2 e^{-|x|^2/4\epsilon}\, dx}{\epsilon^{n/2}\int_\Omega u(x,-\epsilon)^2\, dx} \geq \frac{e^{-n/4\epsilon}}{\epsilon^{n/2}} \quad (4.10)$$

since $|x| \leq \sqrt{n}$ for $x \in \Omega$. We then take the logarithm of (4.10) to get

$$\log \frac{\|u(\cdot,-\epsilon)\|^2_{L^2(\Omega)}}{\|U(\cdot,\tau_0)\|^2} \leq -\frac{n}{2}\log\frac{1}{\epsilon} + \frac{n}{4\epsilon} \lesssim \frac{1}{\epsilon}. \quad (4.11)$$

Using (4.5), (4.9), and (4.11), we obtain

$$\begin{aligned}\log \frac{\|u(\cdot,t)\|^2_{L^2(\Omega)}}{\|U(\cdot,\tau)\|^2} &\lesssim (M_1^2 + M_0)\epsilon + \epsilon^{-1} + (\tau - \log\epsilon^{-1})\epsilon^{-1} \\ &\quad + \epsilon^{1-\beta}M_1(\epsilon^{-1} + \epsilon|a|^2 + 1)^\beta + \epsilon^{1-\alpha}M_0(\epsilon^{-1} + \epsilon|a|^2 + 1)^\alpha \\ &\lesssim M_1^2 + M_0 + \epsilon^{-1} + \tau\epsilon^{-1} + M_1\epsilon^{1-2\beta} + M_0\epsilon^{1-2\alpha},\end{aligned}$$

from where, after an short computation,

$$\log \frac{\|u(\cdot,t)\|^2_{L^2(\Omega)}}{\|U(\cdot,\tau)\|^2} \lesssim M_1^2 + M_0 + M\tau + M_1 M^{2\beta-1} + M_0 M^{2\alpha-1}, \quad (4.12)$$

with $M = M_0^{2/(3-4\alpha)} + M_1^{2/(3-4\beta)}$ as defined in (3.45). Fix $\delta_0 \in (0, 1/2]$. Let $t = -\delta^2$, i.e., $\tau = 2\log(1/\delta)$. If $\delta \in (0, \sqrt{\epsilon})$, then

$$\log \frac{\|u(\cdot,-\delta^2)\|^2_{L^2(\Omega)}}{\|U(\cdot,\log 1/\delta^2)\|^2} \lesssim \gamma(\delta)$$

with $\gamma(\delta)$ defined in (4.4), which equals to the right-hand side of (4.12). Therefore,

$$\|u(\cdot,-\delta)\|^2_{L^2(\Omega)} \lesssim e^{\gamma(\delta)} \int_{\mathbb{R}^n} u(x,-\delta^2)^2 G(x,-\delta^2)\, dx.$$

By Lemma 4.1, we conclude that

$$\|u(\cdot, -\delta^2)\|_{L^2(\Omega)}^2 \lesssim \frac{e^{\gamma(\delta)}}{\delta^n} \int_{B_{\delta_0}} u(x, -\delta^2)^2 \, dx.$$

Our assumption on $\delta$ in (4.2) gives $\gamma(\delta) \lesssim \delta_0^2/\delta^2$ and $\log(1/\delta) \leq \delta_0^2/\delta^2$. So we obtain (4.1) for $t = -\delta^2$. To get (4.1) for $t \in (-\delta^2, 0]$, we simply translate in time. Therefore, we conclude that

$$\|u(\cdot, t)\|_{L^2(\Omega)}^2 \lesssim e^{(n+1)\delta_0^2/\delta^2} \|u(\cdot, t)\|_{L^2(B_{\delta_0})}^2, \qquad t \in [-\delta^2, 0],$$

as desired. $\square$

PROOF OF THEOREM 2.2. We first assume that $u_0$, $v$, and $w$ are smooth as in the proof of Theorem 2.1. Fix any $t \in (0, T)$. For all $\delta_0 \in (0, 1/2]$, by restating Lemma 4.2 to address positive time instead of negative, we claim that

$$\|u(\cdot, s)\|_{L^2(\Omega)}^2 \lesssim e^{P(n, \delta_0, M_0, M_1, p, p)} \|u(\cdot, s)\|_{L^2(B_{\delta_0})}^2, \tag{4.13}$$

for all $s \in [0, \delta^2]$, where $P$ is determined by Lemma 4.2. Choose $t' \in (0, t)$ such that $t - t' < \epsilon$ with $\epsilon$ as in Lemma 4.2. By (4.13), we have

$$\|u(\cdot, t')\|_{L^2(\Omega)}^2 \lesssim e^{P(n, \delta_0, M_0, M_1, q, p)} \|u(\cdot, t')\|_{L^2(B_{\delta_0})}^2.$$

Therefore, (2.10) holds for all $t \in (0, T)$. Thus we have proven Theorem 2.2 under the assumption that $u_0$, $v$, and $w$ are smooth. For the general case, we approximate $u_0$, $v$, and $w$ by smooth functions $u_0^\epsilon$, $v^\epsilon$, and $w^\epsilon$ respectively. We then have $\|u^\epsilon(\cdot, t)\|_{L^2(\Omega)}^2 \to \|u(\cdot, t)\|_{L^2(\Omega)}^2$ and $\|u^\epsilon(\cdot, t)\|_{L^2(B_{\delta_0})}^2 \to \|u(\cdot, t)\|_{L^2(B_{\delta_0})}^2$ as $\epsilon \to 0^+$ by (3.57). Thus (2.10) holds for all $t \in (0, T)$, as desired. $\square$

## 5. The case $\mathbb{R}^n$

In this section, we prove the theorem concerning the case of $\mathbb{R}^n$.

PROOF OF THEOREM 2.3. Note that we now assume the growth condition (2.8) instead of periodicity. The proof is similar to the periodic case with small modifications. With $\epsilon \in (0, T)$ such that $\epsilon \leq 1/4$, we have

$$\frac{\epsilon \int_{\mathbb{R}^n} |\nabla u(x_\epsilon + y, -\epsilon)|^2 G(y, -\epsilon) \, dy}{\int_{\mathbb{R}^n} u(x_\epsilon + y, -\epsilon)^2 G(y, -\epsilon) \, dy} \leq 4K\epsilon q_{\mathrm{D}}(-\epsilon), \tag{5.1}$$

for some $x_\epsilon \in B_2$. This was proven in [**CK**], but since the argument is short, we present it here. To see why (5.1) holds, we assume that

$$\lambda \int_{\mathbb{R}^n} u(x + y, -\epsilon)^2 G(y, -\epsilon) \, dy \leq \int_{\mathbb{R}^n} |\nabla u(x + y, -\epsilon)|^2 G(y, -\epsilon) \, dy, \qquad x \in B_2, \tag{5.2}$$

where $\lambda = 4K q_{\mathrm{D}}(-\epsilon)$. Integrating (5.2) over $B_2$ and proceeding similarly as in (3.1), we have

$$\begin{aligned}
\int_{\mathbb{R}^n} G(y, -\epsilon) \, dy \int_{B_2} u(x + y, -\epsilon)^2 \, dx &\leq \frac{1}{\lambda} \int_{\mathbb{R}^n} G(y, -\epsilon) \, dy \int_{B_2} |\nabla u(x + y, -\epsilon)|^2 \, dx \\
&\leq \frac{1}{\lambda} \|\nabla u(\cdot, -\epsilon)\|^2 \int_{\mathbb{R}^n} G(y, -\epsilon) \, dy \leq \frac{(2\pi)^{n/2}}{\lambda} \|\nabla u(\cdot, -\epsilon)\|^2,
\end{aligned} \tag{5.3}$$

where we continue using the convention (3.12). On the other hand, we have the lower bound for the left hand side of (5.3)

$$\begin{aligned}
\int_{\mathbb{R}^n} G(y, -\epsilon) \, dy \int_{B_2} u(x + y, -\epsilon)^2 \, dx &= \int_{\mathbb{R}^n} G(y, -\epsilon) \, dy \int_{B_2(y)} u(x, -\epsilon)^2 \, dx \\
&= \int_{B_{1/2}} G(y, -\epsilon) \, dy \int_{B_2(y)} u(x, -\epsilon)^2 \, dx + \int_{\mathbb{R}^n \setminus B_{1/2}} G(y, -\epsilon) \, dy \int_{B_2(y)} u(x, -\epsilon)^2 \, dx \\
&\geq \int_{B_{1/2}} G(y, -\epsilon) \, dy \int_{B_2(y)} u(x, -\epsilon)^2 \, dx \geq \int_{B_{1/2}} G(y, -\epsilon) dy \int_{B_{1/2}} G(y, -\epsilon) \, dy \int_{B_1} u(x, -\epsilon)^2 \, dx \\
&\geq \frac{1}{2}(2\pi)^{n/2} \int_{B_1} u(x, -\epsilon)^2 \, dx \geq \frac{1}{2K}(2\pi)^{n/2} \int_{\mathbb{R}^n} u(x, -\epsilon)^2 \, dx,
\end{aligned}$$

where the last inequality holds by the doubling condition (2.8). Thus we obtain

$$\|u(x, -\epsilon)\|^2 \leq \frac{2K}{\lambda}\|\nabla u(\cdot, -\epsilon)\|^2 \leq \frac{1}{2q_{\mathrm{D}}(-\epsilon)}\|\nabla u(\cdot, -\epsilon)\|^2.$$

We have a contradiction since $\|\nabla u(\cdot, -\epsilon)\|^2/\|u(\cdot, -\epsilon)\|^2 \leq q_{\mathrm{D}}(-\epsilon)$. Therefore, (5.2) does not hold for all $x \in B_2$, i.e., there exists $x_\epsilon \in B_2$ such that (5.1) holds true. Theorem 2.3 then follows as in the proofs of Theorems 2.1 and 2.2. □

## 6. $L_t^q L_x^p$ case

The method we used also works for the case where $v \in L_t^{p_2} L_x^{p_1}$ and $w \in L_t^{q_2} L_x^{q_1}$ when $p_2$ and $q_2$ are finite, as shown in the next statement. When limiting $p_2 \to \infty$ and $q_2 \to \infty$, the theorem reduces to the results in Section 2.

THEOREM 6.1. *Let $n \geq 2$. Assume that $u \in L_{x,t}^\infty(\Omega \times I)$ is a solution of (2.1) with $v \in L_t^{p_2} L_x^{p_1}(\Omega \times I)$ and $w \in L_t^{q_2} L_x^{q_1}(\Omega \times I)$ are such that $p_1 > 2n/3$ and $q_1 > 2n$ with*

$$p_2 > \max\left\{\frac{2}{3-4\alpha}, \frac{2}{1-\alpha}\right\} \tag{6.1}$$

*and*

$$q_2 > \max\left\{\frac{2}{3-4\beta}, \frac{2}{1-\beta}\right\},$$

*where $\alpha = n/2p_1$ and $\beta = 1/2 + n/2q_1$. Then for all $(x_0, t_0) \in \mathbb{T}^n \times I$, the vanishing order of $u$ at $(x_0, t_0)$ satisfies*

$$O_{(x_0, t_0)}(u) \lesssim \|v\|_{L_t^{p_2} L_x^{p_1}(\Omega \times I)}^a + \|w\|_{L_t^{q_2} L_x^{q_1}(\Omega \times I)}^b + 1$$

*where*

$$a = \frac{2}{3 - 2/p_2 - 4\alpha} \tag{6.2}$$

*and*

$$b = \frac{2}{3 - 2/q_2 - 4\beta}.$$

PROOF OF THEOREM 6.1. Without loss of generality, we can assume $I = [-1, 0]$ and $(x_0, t_0) = (0, 0)$. For simplicity, we only consider the case $w = 0$ and $v \in L_t^{p_2} L_x^{p_1}$, as for a nonzero $w$ the proof is similar. By a simple change of variables, we obtain

$$\int_{-1}^0 \left(\int |v(x,t)|^{p_1} dx\right)^{p_2/p_1} dt = \int_0^\infty \left(\int |V(y,\tau)|^{p_1} e^{-n\tau/2} dy\right)^{p_2/p_1} e^{-\tau} d\tau = \int_0^\infty K_0(\tau)^{p_2} e^{-\tau} d\tau,$$

where $K_0(\tau) = \|V\|_{L^{p_1}} e^{-n\tau/2p_1}$. Therefore,

$$\int_0^\infty K_0(\tau)^{p_2} e^{-\tau} d\tau \leq M_0^{p_2},$$

where $M_0 = \|v\|_{L_t^{p_2} L_x^{p_1}(\Omega \times I)} + 1$. Note that $K_0(\tau)$ plays the role of $M_0$ in the proof of Theorem 2.1, except that here $K_0$ is a function of $\tau$. As in the proof of Theorem 2.1, we have

$$\begin{aligned}\bar{Q}'(\tau) \lesssim{}& e^{-\tau/2}\epsilon^{-1}\bar{Q}_+^{1/2} + K_0(\tau)^2 e^{(\alpha-2)\tau}\epsilon^{-2\alpha}\bar{Q}_+^\alpha + K_0(\tau)^2 e^{(2\alpha-2)\tau}\bar{Q}_+^{2\alpha} + e^{-\tau}\epsilon^{-2} + e^{-\tau/2}\epsilon^{-1} \\ & + K_0(\tau)^2 e^{-2\tau}\epsilon^{-4\alpha} + K_0(\tau)^2 e^{(\alpha-2)\tau}\epsilon^{-2\alpha} + K_0(\tau)^{2/(1-\alpha)} e^{-2\tau} + K_0(\tau)^2 e^{(2\alpha-2)\tau},\end{aligned} \tag{6.3}$$

where $\alpha = n/2p_1$ and $\epsilon = M_0^{-a}/K$ with $a$ as in (6.2) and $K$ sufficiently large. Denote by $C_2 \geq 1$ the implicit constant in (6.3). Under the condition (6.1), we claim that

$$\bar{Q}_+(\tau) < \frac{C_1 + C_2}{\epsilon},$$

where $C_1 = C_0(q_0 + q_0^{1/2})$ and $C_0 \geq 1$ is sufficiently large. To prove this claim, we use a barrier argument. Assume, contrary to the assertion, that there exists $\tau_1$ such that $\bar{Q}(\tau_1) = (C_1 + C_2)/\epsilon$, and suppose that $\tau_1$ is the first time with

this property. We have

$$\bar{Q}'(\tau) \leq C_2(C_1+C_2)^{1/2}e^{-\tau/2}\epsilon^{-3/2} + C_2K_0(\tau)^2e^{(\alpha-2)\tau}\epsilon^{-3\alpha}(C_1+C_2)^\alpha$$
$$+ C_2K_0(\tau)^2e^{(2\alpha-2)\tau}(C_1+C_2)^{2\alpha}\epsilon^{-2\alpha} + C_2e^{-\tau}\epsilon^{-2} + C_2K_0(\tau)^2e^{-2\tau}\epsilon^{-4\alpha}$$
$$+ C_2K_0(\tau)^2e^{(\alpha-2)\tau}\epsilon^{-2\alpha} + C_2K_0(\tau)^{2/(1-\alpha)}e^{-2\tau} + C_2K_0(\tau)^2e^{(2\alpha-2)\tau}.$$

We integrate (6.3) in $\tau$ from $\tau_0$ to $\tau_1$, where $\tau_0$ is given in (3.10), to obtain

$$\bar{Q}(\tau_1) \leq \bar{Q}(\tau_0) + 2C_2(C_1+C_2)^{1/2}\epsilon^{-1} + C_2(C_1+C_2)^\alpha\epsilon^{-3\alpha}\int_{\tau_0}^{\tau_1} K_0(\tau)^2e^{(\alpha-2)\tau}\,d\tau$$
$$+ C_2(C_1+C_2)^{2\alpha}\epsilon^{-2\alpha}\int_{\tau_0}^{\tau_1} K_0(\tau)^2e^{(2\alpha-2)\tau}d\tau + C_2\epsilon^{-1} + C_2\epsilon^{-4\alpha}\int_{\tau_0}^{\tau_1} K_0(\tau)^2e^{-2\tau}\,d\tau \qquad (6.4)$$
$$+ C_2\epsilon^{-2\alpha}\int_{\tau_0}^{\tau_1} K_0(\tau)^2e^{(\alpha-2)\tau}\,d\tau + C_2\int_{\tau_0}^{\tau_1} K_0(\tau)^{2/(1-\alpha)}e^{-2\tau}\,d\tau + C_2\int_{\tau_0}^{\tau_1} K_0(\tau)^2e^{(2\alpha-2)\tau}\,d\tau.$$

By Hölder's inequality, the third term on the right-hand side satisfies

$$C_2(C_1+C_2)^\alpha\epsilon^{-3\alpha}\int_{\tau_0}^{\tau_1} K_0(\tau)^2e^{(\alpha-2)\tau}\,d\tau$$
$$\leq C_2(C_1+C_2)^\alpha\epsilon^{-3\alpha}\|K_0(\tau)^2e^{-2\tau/p_2}\|_{L^{p_2/2}(0,\infty)}\|e^{(\alpha-2+2/p_2)\tau}\|_{L^{p_2/(p_2-2)}(0,\infty)}$$
$$\leq C_1(C_1+C_2)^\alpha\epsilon^{-3\alpha}\left(\int_0^\infty K_0(\tau)^{p_2}e^{-\tau}\,d\tau\right)^{2/p_2}\epsilon^{2-\alpha-2/p_2} \leq C_1(C_1+C_2)^\alpha\epsilon^{3-4\alpha-2/p_2}M_0^2\epsilon^{-1} \qquad (6.5)$$
$$\leq C_1(C_1+C_2)^\alpha\frac{\epsilon^{-1}}{K^{3-4\alpha-2/p_2}} \leq \frac{C_1+C_2}{100\epsilon},$$

for sufficiently large $C_0$ and $K$, where we used $\epsilon^{3-4\alpha-2/p_2}M_0^2 = 1/K^{3-4\alpha-2/p_2}$ since $a = 2/(3-2/p_2-4\alpha)$ in the second last inequality. In the last inequality, noting that $3-4\alpha-2/p_2 > 0$ by condition (6.1), we have $1/K^{3-4\alpha-2/p_2} < 1/100$ for sufficiently large $K$. We proceed similarly as in (6.5) for every term involving $K_0(\tau)$ on the right-hand side of (6.4) to obtain

$$\bar{Q}(\tau_1) \leq \bar{Q}(\tau_0) + \frac{C_1+C_2}{4\epsilon} < \frac{C_1+C_2}{\epsilon}$$

with $C_0$ and $K$ sufficiently large. Note that the condition $p_2 > 2/(1-\alpha)$ is needed in the tenth term on the right side of (6.4). This is a contradiction, showing that

$$\bar{Q}(\tau) \leq \frac{C_1+C_2}{\epsilon}.$$

Therefore, $\bar{Q}(\tau) \lesssim M_0^a$, as desired. The case $v = 0$ and $w \in L_t^{q_2}L_x^{q_1}$ is similar with $M_0$, $\alpha$, and $K_0$ being replaced by $M_1$, $\beta$, and $K_1$ respectively. The general case, when both $v \in L_t^{p_2}L_x^{p_1}$ and $w \in L_t^{q_2}L_x^{q_1}$ are present, follows analogously. $\square$

## Acknowledgments

The authors were supported in part by the NSF grant DMS-1907992.

ON QUANTITATIVE UNIQUENESS FOR PARABOLIC EQUATIONS 20

DEPARTMENT OF MATHEMATICS, UNIVERSITY OF SOUTHERN CALIFORNIA, LOS ANGELES, CA 90089
*Email address*: kukavica@usc.edu

DEPARTMENT OF MATHEMATICS, UNIVERSITY OF SOUTHERN CALIFORNIA, LOS ANGELES, CA 90089
*Email address*: ntle@usc.edu